\documentclass[a4paper]{article}
\begin{document} 
\newtheorem{prop}{Proposition}[section]
\newtheorem{Def}{Definition}[section] \newtheorem{theorem}{Theorem}[section]
\newtheorem{lemma}{Lemma}[section] \newtheorem{Cor}{Corollary}[section]

\title{\bf Global solutions for the Dirac-Klein-Gordon system in two space 
dimensions}
\author{{\bf Axel Gr\"unrock}\\
Mathematisches Institut\\
Universit\"at Bonn\\
Beringstr. 1\\
53115 Bonn\\
Germany\\
e-mail {\tt gruenroc@math.uni-bonn.de}\\[0.8cm]
{\bf Hartmut Pecher}\\
Fachbereich Mathematik und Naturwissenschaften\\
Bergische Universit\"at Wuppertal\\
Gau{\ss}str.  20\\
42097 Wuppertal\\
Germany\\
e-mail {\tt Hartmut.Pecher@math.uni-wuppertal.de}}
\date{}
\maketitle

\begin{abstract}
The Cauchy problem for the classical Dirac-Klein-Gordon system in two space 
dimensions is globally well-posed for $L^2$ Schr\"odinger data and 
wave data in $H^{\frac{1}{2}} \times H^{-\frac{1}{2}}$ . In the case of smooth 
data there exists  a global smooth (classical) solution. The proof uses 
function spaces of Bourgain type 
based on Besov spaces -- previously applied by Colliander, Kenig and Staffilani 
for generalized Benjamin-Ono equations and also by Bejenaru, Herr, Holmer and 
Tataru for the 2D Zakharov system -- and the null structure of the system 
detected by d'Ancona, Foschi and Selberg, and a refined bilinear Strichartz 
estimate due to Selberg. The global existence proof uses an idea of Colliander, 
Holmer and Tzirakis for the 1D Zakharov system.
\end{abstract}

\renewcommand{\thefootnote}{\fnsymbol{footnote}}
\footnotetext{\hspace{-1.8em}{\it 2000 Mathematics Subject Classification:} 
35Q55, 35L70 
The first author was partially supported by the Deutsche Forschungsgemeinschaft, 
Sonderforschungsbereich 611.}
\normalsize 
\setcounter{section}{0}
\section{Introduction and main results}
Consider the Cauchy problem for the Dirac -- Klein -- Gordon 
equations in two space dimensions 
\begin{eqnarray}
\label{0.1}
i(\partial_t + \alpha \cdot \nabla) \psi + M \beta \psi & = & - \phi \beta \psi 
\\
\label{0.2}
(-\partial_t^2 + \Delta) \phi + m\phi & = & - \langle \beta \psi,\psi \rangle
\end{eqnarray}
with (large) initial data
\begin{equation}
\psi(0)  =  \psi_0 \,  , \, \phi(0)  =  \phi_0 \, , \, \partial_t 
\phi(0) = \phi_1 \, .
\label{0.3}
\end{equation}
Here $\psi$ is a two-spinor field, i.e. $\psi : {\bf R}^{1+2} \to {\bf C}^2$, 
and 
$\phi$ is a real-valued function, i.e. $\phi : {\bf R}^{1+2} \to {\bf R}$ , 
$m,M 
\in {\bf R}$ and $\nabla = (\partial_{x_1} , \partial_{x_2}) $ , $ \alpha \cdot 
\nabla = \alpha^1 \partial_{x_1} + \alpha^2 \partial_{x_2}$ . 
$\alpha^1,\alpha^2, \beta$ are hermitian ($ 2 \times 2$)-matrices satisfying 
$\beta^2 = 
(\alpha^1)^2 = (\alpha^2)^2 = I $ , $ \alpha^j \beta + \beta \alpha^j = 0, $  $ 
\alpha^j \alpha^k + \alpha^k \alpha^j = 2 \delta^{jk} I $ . \\
$\langle \cdot,\cdot \rangle $ denotes the ${\bf C}^2$ - scalar product. A 
particular representation is given by  $\alpha^1 = {0\;\,1 \choose 1\;\,0}$ , 
$\alpha^2 =  {0\,-i \choose i\;\,0}$ , $\beta = {1\;\,0\choose0 -1}$.\\
We consider Cauchy data in Sobolev spaces: $\psi_0 \in H^s $ , $ \phi_0 \in H^r 
$ , $ \phi_1 \in H^{r-1}$ .\\
 The fundamental conservation law is charge conservation $ \| \psi(t) \|_{L^2} 
= 
const $. \\
 In the (1+1)-dimensional case global well-posedness for smooth data was 
already 
established by Chadam \cite{C} and also for much less regular data by 
Bournaveas \cite{B}, Fang \cite{F}, Bournaveas and Gibbeson \cite{BG}, 
Machihara 
\cite{M}, Pecher \cite{P}, Selberg \cite{S1}, Selberg - Tesfahun \cite{ST} and
Tesfahun \cite{T}, the last two authors also for data $\psi_0 \notin L^2$. In 
the 
(2+1)-dimensional and (3+1)-dimensional case no global well-posedness results 
for large data were known so far. In (2+1)-dimensions local well-posedness was 
proven by Bournaveas \cite{B1}, if $s> \frac{1}{4}$ and $r=s+ \frac{1}{2}$, 
which 
was later 
improved by d'Ancona, Foschi and Selberg \cite{AFS1} to the case $ s > - 
\frac{1}{5}$ and $ 
\max(\frac{1}{4}-\frac{s}{2},\frac{1}{4}+\frac{s}{2},s) < r < 
\min(\frac{3}{4}+2s,\frac{3}{4}+\frac{3s}{2},1+s)$. Their proof relied on the 
null structure of the system. This complete null structure was detected by 
d'Ancona, Foschi and Selberg in their earlier paper \cite{AFS}, where it was 
applied to show 
an almost optimal local existence result in (3+1)-dimensions, namely if 
$s=\epsilon$ , $r= \frac{1}{2}+\epsilon$ for any $\epsilon > 0$.
 
We now give the first global well-posedness result for large data in two space 
dimensions. It holds in the case $s=0$ , $r=\frac{1}{2}$, and more generally in 
the case $s \ge 0,$ $r = s +\frac{1}{2}$, where local 
well-posedness was known to be true before already (by d'Ancona, Foschi and 
Selberg \cite{AFS1}). Especially we show the existence of global classical 
solutions for smooth data. It is necessary to refine the local existence 
result by 
replacing Bourgain spaces $X^{s,b}_{\pm}$ and $X^{r,b}_{\pm}$ for $ b > 
\frac{1}{2}$ constructed from Sobolev spaces by their analogue constructed from 
Besov spaces with respect to time, especially $X_{\pm}^{s,\frac{1}{3},1}$ and 
$X_{\pm}^{r,\frac{1}{3},1}$ (see the definition below). Spaces of this type were 
already successfully used to give a local well-posedness result for the 2D - 
Zakharov system by Bejenaru, Herr, Holmer and Tataru \cite{BHHT} and 
Colliander, 
Kenig and Staffilani for generalized Benjamin-Ono equations \cite{CKS}. The 
precise bound for the 
existence time then can be combined with the charge conservation to show global 
well-posedness for our 2D Dirac-Klein-Gordon system. A similar procedure was 
already used by Colliander, Holmer and Tzirakis for the one-dimensional 
Zakharov 
system \cite{CHT}. It turns out that the choice of the regularity parameters 
$s$ 
and $r$ in 
our case just allows to estimate both nonlinearities in a unified way. What one 
also needs are of course the Strichartz estimates for the wave equation, here 
also the Besov space version to avoid the endpoint Strichartz estimate in 2D. 
The Strichartz estimates however are not sufficient for a particularly delicate 
case where it is essential to use a bilinear refinement which was detected by 
Selberg \cite{S} and can also be found in Foschi-Klainerman \cite{FK}. This 
version was already 
used by d'Ancona, Foschi and Selberg \cite{AFS1} in their local well-posedness 
result.

We use the following function spaces. Let $\, \widehat{} \,$ denote the Fourier 
transform with respect to space or time and $\, \tilde{} \,$ the Fourier 
transform with respect to space and time simultaneously. Let $\varphi \in 
C_0^{\infty}({\bf R^n})$ 
be a nonnegative function with $ supp \, \varphi \subset \{ 1/2 \le |\xi| \le 2 
\} $ and $\varphi(\xi) > 0 $ , if $ \frac{1}{\sqrt{2}} \le |\xi| \le \sqrt{2} 
$. 
Setting $\widehat{\rho}_k(\xi) := \varphi(2^{-k} \xi)$ $(k=1,2,...)$ , 
$\widehat{\varphi}_k(\xi) := 
\frac{\widehat{\rho}_k(\xi)}{\sum_{j=-\infty}^{+\infty} \varphi(2^{-j}\xi)}$ 
$(k=1,2,...)$ and $ \widehat{\varphi}_0(\xi):= 1 - \sum_{k=1}^{\infty} 
\widehat{\varphi}_k(\xi)$ we have $supp \, \widehat{\varphi}_k \subset 
\{2^{k-1} 
\le |\xi| \le 2^{k+1}\}$ , $ supp \, \widehat{\varphi}_0 \subset \{ |\xi| \le 2 
\}$ and $\sum_{k=0}^{\infty} \widehat{\varphi}_k = 1 $. The Besov spaces are 
defined for $s \in {\bf R}$ , $1 \le p,q \le \infty$ as follows:
$$ B^s_{p,q} = \{ f \in {\cal S}' \, , \, \|f\|_{B^s_{p,q}} < \infty \} \, , $$ 
where
$$ \|f\|_{B^s_{p,q}} = \left( \sum_{k=0}^{\infty} (2^{sk} \|\varphi_k \ast f 
\|_{L^p})^q \right)^{\frac{1}{q}}  {\mbox if} \,q < \infty \, , $$
$$ \|f\|_{B^s_{p,\infty}} = \sup_{k \ge 0} 2^{sk} \|\varphi_k \ast f \|_{L^p} 
$$
(cf. e.g. Triebel \cite{Tr}, Section 2.3.1).\\
Similarly the homogeneous Besov spaces are defined as the set of those $f\in 
{\cal S}'$, for which $\|f\|_{\dot{B}^s_{p,q}}$ is finite, where
$ \|f\|_{\dot{B}^s_{p,q}} = (\sum_{k=-\infty}^{+\infty} (2^{sk} \|\varphi_k '
 \ast f\|_{L^p})^q)^{\frac{1}{q}} $ with the usual modification for $q=\infty$ 
and $ \widehat{\varphi}_k'(\xi) := 
\frac{\widehat{\rho}_k(\xi)}{\sum_{j=-\infty}^{+\infty} \varphi(2^{-j}\xi)}$ 
for 
$k\in {\bf Z}$. 
We also need the following Bourgain type spaces. The standard spaces belonging 
to 
the half waves are defined by the completion of ${\cal S}({\bf R} \times {\bf 
R^2})$ with respect to
$$ \|f\|_{X^{s,b}_{\pm}} = \|U_{\pm}(-t)f\|_{H^b_t H^s_x} = \| \langle \xi 
\rangle^s \langle \tau \pm |\xi| \rangle^b \tilde{f}(\tau,\xi) \|_{L^2} \, $$
where 
$$ U_{\pm}(t):=e^{\mp it|D|} \quad {\mbox {and}} \quad \|g\|_{H^b_t H^s_x} = \| 
\langle 
\xi \rangle^s \langle \tau \rangle^b \tilde{g}(\xi,\tau)\|_{L^2_{\xi,\tau}} \, 
. 
$$
We also define $X^{s,b,q}_{\pm}$ as the space of all $u \in {\cal S}'({\bf R} 
\times {\bf R^2})$, where the following norms are finite:
$$ \|f\|_{X^{s,b,q}_{\pm}} = \|U_{\pm}(-t)f\|_{B^b_{2,q} H^s_x} = \left( 
\sum_{k=0}^{\infty} 2^ {qbk} \| \langle \xi \rangle^s \widehat{\varphi}_k(\tau 
\pm |\xi|) \tilde{f}(\tau,\xi) \|_{L^2_{\tau,\xi}}^q \right)^{\frac{1}{q}} $$
for $ 1 \le q < \infty$ , where
$$ \|g\|_{B^b_{2,q}H^s_x}= \left(\sum_{k=0}^{\infty} 2^{qbk} \|\langle \xi 
\rangle^s \widehat{\varphi}_k(\tau) \tilde{f}(\tau,\xi)\|_{L^2_{\tau,\xi}}^q 
\right)^{\frac{1}{q}}$$
and
$$ \|f\|_{X^{s,b,\infty}_{\pm}} = \|U_{\pm}(-t) f\|_{B^b_{2,\infty} H^s_x} = 
\sup_{k \ge 0} 2^{bk} \| \langle \xi \rangle^s \widehat{\varphi}_k(\tau \pm 
|\xi|) \tilde{f}(\tau,\xi)\|_{L^2_{\tau,\xi}}$$
for $ q = \infty$ , where
$$\|g\|_{B^b_{2,\infty} H^s_x} = \sup_{k \ge 0} 2^{bk} \| \langle \xi \rangle^s 
\widehat{\varphi}_k(\tau) \tilde{g}(\tau,\xi)\|_{L^2_{\tau,\xi}}  \, .$$
Note that $U_{\pm}(t) = e^{\mp it(-\Delta +1)^{1/2}}$ would lead to equivalent 
norms.

Spaces of type $X^{s,b,q}$ with various phase functions $\phi(\xi)$ instead of 
$\pm |\xi|$ have been used
in the literature before, for example by Colliander, Kenig, and Staffilani in 
their work on dispersion
generalized Benjamin-Ono equations \cite{CKS}. As was observed in \cite[proof 
of Lemma 5.1]{CKS}, they can
be obtained by real interpolation from the standard $X^{s,b}$-spaces. In fact, 
by \cite[Theorem 5.6.1]{BL}
one has for $s \in {\bf R}$, $1 \le q \le \infty$, $b_0 \neq b_1$ and 
$b=(1-\theta)b_0 + \theta b_1$, $0<\theta<1$, that
$$ (X^{s,b_0},X^{s,b_1})_{\theta,q}=X^{s,b,q}.$$
Using the duality Theorem \cite[Theorem 3.7.1]{BL} we see that for $1 \le q < 
\infty$
$$ (\overline{X}^{s,b,q})' = X^{-s,-b,q'},$$
where $\overline{X}$ denotes the space of complex conjugates of elements of $X$ 
with norm $\|f\|_{\overline{X}}=\|\overline{f}\|_X$. In the proof of the 
crucial bilinear estimates for local 
well-posedness we will repeatedly make use of complex interpolation. To justify 
this we use the corresponding 
theorem on interpolation of spaces of vector valued sequences \cite[Theorem 
5.6.3]{BL} and take into account
the considerations in \cite[Section 6.4]{BL} to see that
$$ (X^{s_0,b_0, q_0},X^{s_1,b_1,q_1})_{[\theta]}=X^{s,b,q},$$
whenever $0<\theta  <1$, $s=(1-\theta)s_0 + \theta s_1$, $b=(1-\theta)b_0 + 
\theta b_1$, and $1 \le q_0,q_1 \le \infty$ as well as $\frac1q = 
\frac{1-\theta}{q_0}+\frac{\theta}{q_1}$.\\
The preceeding remarks on duality and interpolation are completely independent 
of the specific phase function.\\
For 
$B \subset {\cal S}^{'}({\bf R} \times {\bf R^2})$ we denote by $B(T)$ 
the space of restrictions of distributions in $B$ to the set $(0,T) \times {\bf 
R^2}$  with induced norm.

We use the Strichartz estimates for the homogeneous wave equation in ${\bf R^n} 
\times {\bf R}$, which can be found e.g. in Ginibre-Velo \cite{GV}, Prop. 2.1.
\begin{prop}
Let $\gamma(r) = (n-1)(\frac{1}{2}-\frac{1}{r})$ , $ \delta(r) = 
n(\frac{1}{2}-\frac{1}{r}) $ , $ n \ge 2 $. Let $\rho,\mu \in {\bf R}$, $2 \le 
q,r \le \infty$ satisfy $0 \le \frac{2}{q} \le \min(\gamma(r),1) $ , $ 
(\frac{2}{q},\gamma(r)) \neq (1,1)$, $ \rho + \delta(r) - \frac{1}{q} = \mu $. 
Then
$$ \|e^{\pm it|D|}u_0\|_{L^q({\bf R}, \dot{B}^{\rho}_{r,2}({\bf R}^n))} \le c 
\|u_0\|_{\dot{H}^{\mu}({\bf R}^n)} \, . $$
The same holds with $\dot{B}^{\rho}_{r,2}$ replaced by $\dot{H}^{\rho,r}$ under 
the additional assumption $r < \infty$. 
\end{prop}
The following consequence of estimates of Strichartz type is important for our 
considerations.
\begin{prop}
\label{Prop. 1.1}
Let $Y \subset {\cal S}^{'}({\bf R} \times {\bf R}^n)$ be a set of functions of 
space and time with the property that
$$ \|hf\|_Y \le c \|h\|_{L^{\infty}_t} \|f\|_Y $$
for all $h\in L^{\infty}_t$ and $f \in Y$. Assume moreover the (Strichartz 
type) 
estimate
$$ \|U_{\pm}(t)u_0\|_Y \le c \|u_0\|_{H^{\mu}} \,  , $$
where $ U_{\pm}(t) = e^{\mp it|D|}$. Then the following estimate holds:
$$ \|f\|_Y \le c \|f\|_{X^{\mu,\frac{1}{2},1}_{\pm}} \, . $$
\end{prop}
{\bf Proof:}
We combine Lemma 2.3 in \cite{GTV} with the proof of the embedding $ 
B^{1/2}_{2,1}({\bf R})$ $\subset C^0({\bf R})$. Let $\psi$ be a 
$C_0^{\infty}({\bf 
R})$ - function with $\psi(\tau) = 1 $ for $ 1/2 \le |\tau| \le 2$ and
$\widehat{\psi}_k(\tau) := \psi(2^{-k} \tau)$ , so that $\widehat{\psi}_k 
(\tau) = 1$ 
for $ 2^{k-1} \le |\tau| \le 2^{k+1}$. Furthermore we define $\widehat{\psi}_0 
\in C^{\infty}_0$ such that $\widehat{\psi}_0(\tau) = 1 $ for $|\tau| \le 2$.
The functions $\varphi_k$ are those which appear in the definition of the Besov 
norms (here in the 1-dimensional case). We thus have the property that 
$\widehat{\psi_k}(\tau) = 1 $ for $ \tau \in supp \, \widehat{\varphi}_k$ 
$(k=0,1,2,...)$. We start from
$$ f = \int e^{it\tau} U_{\pm}(t) ( {\cal F}_t U_{\pm}(- \cdot) f)(\tau) d \tau 
\, . $$
Then we have with $ h = e^{it\tau} $ (for fixed $\tau$):
\begin{eqnarray*}
\|f\|_Y & \le & \int \|U_{\pm}(t)({\cal F}_t U_{\pm}(- \cdot) f)(\tau) \|_Y 
d\tau 
 \le  c \int \|{\cal F}_t U_{\pm}(- \cdot)f\|_{H^{\mu}} d\tau \\
& \le & c \sum_{k=0}^{\infty} \|{\cal F}_t(U_{\pm}(-\cdot)f) 
\widehat{\varphi}_k\|_{L^1_\tau(H^{\mu}_x)} 
 = c \sum_{k=0}^{\infty} \|{\cal F}_t(U_{\pm}(-\cdot)f) \widehat{\varphi}_k 
\widehat{\psi}_k\|_{L^1_\tau(H^{\mu}_x)} \\
& \le & c \sum_{k=0}^{\infty} \|{\cal F}_t(U_{\pm}(-\cdot) f \ast_t 
\varphi_k)(\tau) \|_{L^2_{\tau}(H^{\mu}_x)} \|\widehat{\psi}_k\|_{L^2_{\tau}} 
\\
& \le & c \sum_{k=0}^{\infty} 2^{k/2} \|U_{\pm}(- \cdot)f \ast_t 
\varphi_k\|_{L^2_t(H^{\mu}_x)} \\
& = & c \|U_{\pm}(-t)f\|_{B^{1/2}_{2,1} H^{\mu}_x} 
 =  c \|f\|_{X^{\mu,\frac{1}{2},1}_{\pm}} \, . 
\end{eqnarray*}
where we used $ \|\widehat{\psi}_k\|_{L^2} = 2^{k/2} \|\psi\|_{L^2} $ , 
$(k=1,2,...)$ .

Similarly one can prove a bilinear version:
\begin{prop}
\label{Prop. 1.2}
Let $Y$ be as in Prop. \ref{Prop. 1.1}. Assume
$$ \|U_{\pm 1}(t)u_0 U_{\pm 2}u_1\|_Y \le c \|u_0\|_{H^{\mu}} 
\|u_1\|_{H^{\lambda}} \, . $$
Then
$$ \|f_0 f_1\|_Y \le c \|f_0\|_{X^{\mu,\frac{1}{2},1}_{\pm 1}} 
\|f_1\|_{X^{\lambda,\frac{1}{2},1}_{\pm 2}}
 \, . $$
Here $\pm_1 $ and $\pm_2$ denote independent signs.
\end{prop}
The main result reads as follows:
\begin{theorem}
\label{Theorem}
The Cauchy problem for the Dirac - Klein - Gordon system (\ref{0.1}), 
(\ref{0.2}), (\ref{0.3}) is globally well-posed for data
$$ \psi_0 \in L^2({\bf R^2}), \phi_0 \in H^{1/2}({\bf R^2}), \phi_1 \in 
H^{-1/2}({\bf R^2}) \, . $$
More precisely there exists a unique global solution $(\psi,\phi)$ such that for 
all $ T>0$
\begin{eqnarray*} \psi \in X^{0,\frac{1}{3},1}_+(T) + X^{0,\frac{1}{3},1}_-(T) 
\, , \, \phi \in  
X^{\frac{1}{2},\frac{1}{3},1}_+(T) + X^{\frac{1}{2},\frac{1}{3},1}_-(T) \, , \\ 
\partial_t \phi \in  X^{-\frac{1}{2},\frac{1}{3},1}_+(T) + 
X^{-\frac{1}{2},\frac{1}{3},1}_- (T)\, . \hspace{2cm}
\end{eqnarray*}
This solution has the property
$$ \psi \in C^0({\bf R^+},L^2({\bf R}^2)) \, , \, \phi \in C^0({\bf 
R^+},H^{\frac{1}{2}}({\bf R}^2)) \, , \, \partial_t \phi \in C^0({\bf 
R^+},H^{-\frac{1}{2}}({\bf R}^2)) \, . $$
\end{theorem}
For more regular data we also get the following result.
\begin{theorem} \label{Theorem'}
Let $s$ be an arbitrary nonnegative number. If
 $$\psi_0 \in H^s({\bf R}^2)\, , \, \phi_0 \in H^{s+\frac{1}{2}}({\bf R}^2) \, 
, \, \phi_1 \in H^{s-\frac{1}{2}}({\bf R}^2)\, ,$$
  the global solution of Theorem \ref{Theorem} has the properties:
For every $T>0$
\begin{eqnarray*} 
\psi \in X^{s,\frac{1}{3},1}_+(T) + X^{s,\frac{1}{3},1}_-(T) \, , \, 
 \phi \in X^{s+\frac{1}{2},\frac{1}{3},1}_+(T) + 
X^{s+\frac{1}{2},\frac{1}{3},1}_-(T) 
\, , \\ \partial_t \phi \in X^{s-\frac{1}{2},\frac{1}{3},1}_+(T) + 
X^{s-\frac{1}{2},\frac{1}{3},1}_-(T) \, \hspace{2cm}
\end{eqnarray*}
and

$$ \psi \in C^0({\bf R^+},H^s({\bf R}^2)) \, , \phi \in C^0({\bf 
R^+},H^{s+\frac{1}{2}}({\bf R}^2)) \, , \partial_t \phi \in C^0({\bf 
R^+},H^{s-\frac{1}{2}}({\bf R}^2))\, . $$

If $ s > \frac{5}{2}$ , this is a classical solution, i.e. 
$$ \psi \in C^1({\bf R^+} \times {\bf R}^2) \quad , \quad \phi \in C^2({\bf 
R^+} \times {\bf R}^2)\, . $$
 \end{theorem}
  
\section{Proof of the Theorems}
It is possible to simplify the system (\ref{0.1}),(\ref{0.2}),(\ref{0.3}) by 
considering the projections onto the one-dimensional eigenspaces of the 
operator 
$-i \alpha \cdot \nabla$ belonging to the eigenvalues $ \pm |\xi|$. These 
projections are given by $\Pi_{\pm}(D)$, where  $ D = 
\frac{\nabla}{i} $ and $\Pi_{\pm}(\xi) = \frac{1}{2}(I 
\pm \frac{\xi}{|\xi|} \cdot \alpha) $. Then $ 
-i\alpha \cdot \nabla = |D| \Pi_+(D) - |D| \Pi_-(D) $ and $ \Pi_{\pm}(\xi) \beta
= \beta \Pi_{\mp}(\xi) $. Defining $ \psi_{\pm} := \Pi_{\pm}(D) \psi$ and 
splitting the function $\phi$ into the sum $\phi = \frac{1}{2}(\phi_+ + 
\phi_-)$, where $\phi_{\pm} := \phi \pm iA^{-1/2} \partial_t \phi $ , $ A:= 
-\Delta+1$ , the Dirac - Klein - Gordon system can be rewritten as
\begin{eqnarray}
\label{*}
(-i \partial_t \pm |D|)\psi_{\pm} & = & -M\beta \psi_{\mp} + \Pi_{\pm}(\phi 
\beta (\psi_+ + \psi_-)) \\
\nonumber
(i\partial_t \mp A^{1/2})\phi_{\pm} & = &\mp A^{-1/2} \langle \beta (\psi_+ + 
\psi_-), \psi+ + \psi_- \rangle \mp A^{-1/2} (m+1)(\phi_+ + \phi_-) . \\
\label {**}
\end{eqnarray}
The initial conditions are transformed into
\begin{equation}
\label{***}
\psi_{\pm}(0) = \Pi_{\pm}(D)\psi_0 \, , \,  \phi_{\pm}(0) = \phi_0 \pm i 
A^{-1/2} \phi_1
\end{equation}
In the following we consider the system of integral equations belonging to the 
Cauchy problem (\ref{*}),(\ref{**}),(\ref{***}):
\begin{eqnarray}
\nonumber
\psi_{\pm}(t) & \hspace{-0.3cm} = & \hspace{-0.3cm} e^{\mp it|D|} \psi_{\pm}(0) 
-i \int_0^t e^{\mp i(t-s)|D|} 
\Pi_{\pm}(D)(\frac{1}{2}(\phi_+(s)+\phi_-(s))\beta(\Pi_+(D) \psi_+(s) \\
\label{1.3}
& & 
\hspace{2.5cm} + \Pi_-(D)\psi_-(s)) ds + iM \int_0^t e^{\mp i(t-s)|D|} \beta 
\psi{\mp}(s) ds \\
\nonumber
 \phi_{\pm}(t)& \hspace{-0.3cm}= & \hspace{-0.3cm}e^{\mp itA^{1/2}} 
\phi_{\pm}(0) \pm i \int_0^t e^{\mp 
i(t-s)A^{1/2}} A^{-1/2} \langle \beta(\Pi_+(D) \psi_+(s)+\Pi_-(D)\psi_-(s), \\
 \label{1.4}
& & \hspace{2.5cm} \Pi_+(D)\psi_+(s) + \Pi_-(D)\psi_-(s) \rangle ds \\
\nonumber
& & \hspace{2.5cm} \pm i(m+1)\int_0^t e^{\mp i(t-s)A^{1/2}} A^{-1/2}(\phi_+(s) 
+ 
\phi_-(s)) ds
\end{eqnarray}
We remark that any solution of this system automatically fulfills 
$\Pi_{\pm}(D)\psi_{\pm} = \psi_{\pm}$, because applying $\Pi_{\pm}(D)$ to the 
right hand side of (\ref{1.3}) gives $\Pi_{\pm}(D)\psi_{\pm}(0) = 
\psi_{\pm}(0)$ 
and the integral terms also remain unchanged, because $\Pi_{\pm}(D)^2 = 
\Pi_{\pm}(D)$ and $\Pi_{\pm}(D) \beta \psi_{\mp}(s) = \beta 
\Pi_{\mp}(D)\psi_{\mp}(s) = \beta\psi_{\mp}(s)$. Thus $\Pi_{\pm}(D)\psi_{\pm}$ 
can be 
replaced by $\psi_{\pm}$, thus the system of integral equations reduces exactly 
to the one belonging to our Cauchy problem (\ref{*}),(\ref{**}),(\ref{***}).
 
In order to construct solutions to this system of integral equations we use the 
following facts for the linear problem which are independent of the specific 
phase function.
The following Proposition is closely related to the 
exposition in \cite{BHHT}[Section 5], where slightly different function spaces 
are considered. For the moment let $\psi$ denote
a smooth time cut-off function and set $\psi_T(t)=\psi(\frac{t}{T})$, 
where $0< T \le 1$. The solution
of the inhomogeneous linear equation
$$\partial _t v -i \phi (D) v=F\hspace{1cm}v(0)=0$$
is denoted by $U_{*R} F$, defined by
$$U_{*R} F(t)= \int_0^tU(t-t')F(t')dt',$$
where $U(t)u_0 = e^{it\phi(D)}u_0$ solves the corresponding homogeneous 
equation with initial datum $u_0$
(cf. \cite[Section 2]{GTV}).

\begin{prop}\label{linest}
Let $0<T\le 1$, $-\frac12 < b' < 0 < b \le \frac12$, $s \in {\bf R}$, $u_0 \in 
H^s$ and $F \in L_t^1(I,H^s)$ for a time
interval $I \subset {\bf R}$. Then
\begin{itemize}
 \item[i)] $\|\psi_T U u_0\|_{X^{s,\frac12,1}} \le c \|u_0\|_{H^s}$,
 \item[ii)] $\|\psi_T U _{*R} F\|_{X^{s,\frac12,1}} \le c T^{\frac12 +b'} 
\|F\|_{X^{s,b',\infty}}$,
 \item[iii)] $\|\psi_T u\|_{X^{s,b,1}} \le c T^{\frac12 -b} 
\|u\|_{X^{s,\frac12,1}}$.
\end{itemize}
Moreover $X^{s,\frac12,1} \subset C^0({\bf R}, H^s)$ with a continuous 
embedding.
\end{prop}
{\bf Proof:}
 Without loss of generality we may assume $s=0$. Then we consider the scaling 
transformations $S_T$
and $S^T$ defined by
$$S_Tf(t,x)= f(\frac{t}{T},x), \quad S^Tf(t,x)=Tf(Tt,x),$$
which are formally adjoint to each other with respect to the inner product in 
$L^2_{tx}$ (or merely in $L^2_t$,
if $f$ does not depend on $x$). An elementary calculation shows that for $b >0$
\begin{equation}
 \label{scale}
\|S_T f\|_{B^b_{2,q}L^2_x} \le c T^{\frac12 -b}\| f\|_{B^b_{2,q}L^2_x},
\end{equation}
which is still true without the additional $L^2_x$-part of the norm. For $b = 
\frac12$ and $q=1$ we especially
obtain i), when replacing $f$ by $\psi u_0$. To see ii), we write $Kf(t) = 
\int_0^t f(t')dt'$. Then  $\psi_T Kf(Tt) = \psi K S^Tf(t)$ (cf. \cite[p. 
20]{BHHT}),
hence $\psi_T Kf=S_T(\psi K S^Tf)$ and thus
\begin{eqnarray*}
 \|\psi_T Kf\|_{B^{\frac12}_{2,1}L^2_x} &=& \|S_T(\psi  
KS^Tf)\|_{B^{\frac12}_{2,1}L^2_x} \le c\|\psi  KS^Tf\|_{B^{\frac12}_{2,1}L^2_x} 
\\
& \le c & \|\psi  KS^Tf\|_{H^{\frac12+}L^2_x}
 \le c  \|  KS^Tf\|_{H^{\frac12+}L^2_x} \\
&\le c & \|  S^Tf\|_{H^{-\frac12+}L^2_x} \le c \|  
S^Tf\|_{B^{b'}_{2,\infty}L^2_x},
\end{eqnarray*}
where $b'> -\frac12$. Here we used (\ref{scale}), $H^{\frac12+} \subset 
B^{\frac12}_{2,1}$, the fact that $H^{\frac12+}$ forms an algebra,
Lemma 2.1 from \cite{GTV}, and $B^{b'}_{2,\infty}\subset H^{b'-}$. Dualizing 
(\ref{scale}) we see that the latter is bounded
by $T^{\frac12 +b'} \|f\|_{B^{b'}_{2,\infty}L^2_x}$, which gives ii), when 
replacing $f$ by $U(- \cdot)F$. Part iii) is a consequence of (\ref{scale}) and 
the multiplication
law for $1$-Besov-spaces below. The additional statement follows from the 
well-known embedding $B^{\frac12}_{2,1}\subset C^0$.

\begin{lemma}[One-dimensional Besov-multiplication-law]
\label{multlaw}
 For $0<b \le \frac12$ we have
$$ \|\psi u\|_{B^b_{2,1}L^2_x} \le c \|\psi \|_{B^b_{2,1}}\| 
u\|_{B^{\frac12}_{2,1}L^2_x}.$$
\end{lemma}

{\bf Proof:}
 Let $P_ku = \varphi_k * u$, where $\varphi_k$ are the defining functions of 
the Besov spaces, and $\widetilde{P}_k = P_{k-1}+P_k+P_{k+1}$. Then
\begin{eqnarray*}
 \|\psi u\|_{B^b_{2,1}L^2_x} & = &\ \sum_{l \ge 0} 2^{lb}\|P_l (\psi 
u)\|_{L^2_{tx}} \le \sum_{k,l \ge 0} 2^{lb}\|P_l ((P_k\psi) u)\|_{L^2_{tx}}\\
& \le &c  (\sum_{l \le k+2} 2^{lb}\|(P_k\psi) u\|_{L^2_{tx}}+\sum_{l \ge k+3} 
2^{lb}\|(P_k\psi)(\widetilde{P}_l u)\|_{L^2_{tx}})=:\hspace{-2pt} \sum_1 + 
\sum_2
\end{eqnarray*}
with
$$\sum_1 \le c \sum_{k \ge 0} 
2^{kb}\|P_k\psi\|_{L^2_{t}}\|u\|_{L_t^{\infty}L^2_x}\le c \|\psi 
\|_{B^b_{2,1}}\| u\|_{B^{\frac12}_{2,1}L^2_x},$$
since $B^{\frac12}_{2,1} \subset L^{\infty}$. To estimate $\sum_2$ we choose 
$\frac1p = \frac12 - b$, $\frac1q = b$ so that
\begin{eqnarray*}
 \sum_2 & \le c & \sum_{k \ge 0} \|P_k\psi\|_{L^p_{t}}\sum_{l \ge 0}2^{lb}\|P_l 
u\|_{L_t^qL^2_x}\\
& \le c & \sum_{k \ge 0} \|\varphi_k\|_{L_t^{q'}}\|\widetilde{P}_k 
\psi\|_{L^2_{t}}\sum_{l \ge 0}2^{lb}\|\varphi_l\|_{L_t^{p'}}\|\widetilde{P}_l 
u\|_{L^2_{tx}},
\end{eqnarray*}
where we used Young's inequality. Since $\|\varphi_k\|_{L_t^{q'}} \sim 
2^{\frac{k}{q}}$, we obtain the desired bound.

\begin{prop}
\label{Prop. 2.2'}
For $ 0 \le b' < 1/2$ and $0 < T \le  1$ we have
$$ \|f\|_{L^2({\bf R^2}\times [0,T])} \le c T^{b'} \|f\|_{X^{0,b'}} $$
and
$$ \|f\|_{X^{0,-b'}} \le c T^{b'} \|f\|_{L^2({\bf R^2}\times [0,T])} \, . 
$$
\end{prop}
{\bf Proof:} 
 By the embedding $ H^{b'} \subset L^{\frac{2}{1-2b'}}$  $(0 \le b' 
< 1/2)$
we get
$$ \|\psi_T g \|_{L^2[0,T]} \le \|\psi_T\|_{L^{\frac{1}{b'}}} 
\|g\|_{L^{\frac{2}{1-2b'}}} \le c T^{b'}\|\psi\|_{L^{\frac{1}{b'}}} 
\|g\|_{H^{b'}} \, . $$
From this we get:
\begin{eqnarray*}
\|\psi_T f \|_{L^2_{xt}} & = & \|U(- \cdot) \psi_T f\|_{L^2_{xt}} = 
\|\psi_T U(-\cdot)f\|_{L^2_{xt}} \\
& \le & c T^{b'} \|U(- \cdot) f\|_{H^{b'}_t L^2_x} = c T^{b'} 
\|f\|_{X^{0,{b'}}} \, , 
\end{eqnarray*}
The second claim follows by duality.

Concerning the nonlinearities we shall prove the following estimates in Chapter 
3 below. Here and in the sequel the letter $\psi$ is used again to denote the 
spinor field.
\begin{prop}
\label{Prop. 2.3}
The following estimates are true:
\begin{equation}
\label{1.1}
\|\langle \beta \Pi_{\pm 1}(D)\psi,\Pi_{\pm 
2}(D)\psi^{'}\rangle\|_{X^{-\frac{1}{2},-\frac{1}{3},\infty}_{\pm 3}} \le c 
\|\psi\|_{X^{0,\frac{1}{3},1}_{\pm 1}} \|\psi^{'}\|_{X^{0,\frac{1}{3},1}_{\pm 
2}} 
\end{equation}
and 
\begin{equation}
\label{1.2}
\|\Pi_{\pm 2}(D)(\phi \beta \Pi_{\pm 
1}(D)\psi)\|_{X^{0,-\frac{1}{3},\infty}_{\pm 2}} \le c 
\|\phi\|_{X^{\frac{1}{2},\frac{1}{3},1}_{\pm 3}} 
\|\psi\|_{X^{0,\frac{1}{3},1}_{\pm 1}} \, .
\end{equation}
Here and in the following $\pm_1, \pm_2, \pm_3$ denote independent signs.
\end{prop}

The following local existence result now is a consequence of these estimates.
\begin{prop}
\label{Prop. 2.4}
Let $ \psi_{\pm}(0) \in L^2({\bf R^2}) $ , $ \phi_{\pm}(0) \in 
H^{\frac{1}{2}}({\bf R^2})$ . Then there exists $ 1 \ge T >0$ such that the 
system of integral equations (\ref{1.3}),(\ref{1.4}) has a unique solution
$$ \psi_{\pm} \in X^{0,\frac{1}{3},1}_{\pm}(T) \, , \, \phi_{\pm} \in 
X^{\frac{1}{2},\frac{1}{3},1}_{\pm}(T) \, . $$
This solution has the following properties:
$$ \psi_{\pm}\ \in C^0([0,T],L^2({\bf R^2})) \, , \, \phi_{\pm} \in C^0([0,T], 
H^{\frac{1}{2}}({\bf R}^2)) \, . $$
$\phi_{\pm}$ fulfills
\begin{eqnarray} \label{*****}
& &\|\phi_+(t)\|_{H^{\frac{1}{2}}_x} + \|\phi_-(t)\|_{H^{\frac{1}{2}}_x}\\ 
\nonumber
& \le & (\|\phi_+(0)\|_{H^{\frac{1}{2}}} + \|\phi_-(0)\|_{H^{\frac{1}{2}}}) + 
cT^{\frac{1}{2}} (\|\psi_+(0)\|_{L^2}^2 + \|\psi_-(0)\|_{L^2}^2) + c_0 
T^{\frac{1}{2}} \, ,
\end{eqnarray}
for $0 \le t \le T$, where $c_0$ is a fixed constant.
$T$ can be chosen such that
\begin{equation}
\label{1.5}
T^{\frac{1}{2}}(\|\psi_+(0)\|_{L^2} + \|\psi_-(0)\|_{L^2}) \le c \, ,
\end{equation}
\begin{equation}
\label{1.6}  
T^{\frac{1}{2}}(\|\phi_+(0)\|_{H^{\frac{1}{2}}} + 
\|\phi_-(0)\|_{H^{\frac{1}{2}}}) \le c \, ,
\end{equation}
\begin{equation}
\label{1.7}
T^{\frac{1}{2}}(\|\psi_+(0)\|_{L^2}^2 + \|\psi_-(0)\|_{L^2}^2) \le 
c(\|\phi_+(0)\|_{H^{\frac{1}{2}}} + \|\phi_-(0)\|_{H^{\frac{1}{2}}})\, .
\end{equation}
In addition, if $T$ fulfills only (\ref{1.5}) and (\ref{1.6}) we get the same 
result except estimate (\ref{*****}).
\end{prop}
{\bf Proof:}
Consider the transformation mapping the left hand side of our integral 
equations 
(\ref{1.3}),(\ref{1.4}) into the right hand sides. We construct a fixed point 
of 
it by the contraction mapping principle in the following set
\begin{eqnarray*}
M_T & := & \{ \psi_{\pm} \in X^{0,\frac{1}{3},1}_{\pm}(T) \, , \, \phi_{\pm} \in 
X^{\frac{1}{2},\frac{1}{3},1}_{\pm}(T) : \\
& & \|\psi_+\|_{X^{0,\frac{1}{3},1}_+} + \|\psi_-\|_{X^{0,\frac{1}{3},1}_-} \le 
c T^{\frac{1}{6}} (\|\psi_+(0)\|_{L^2} + \|\psi_-(0)\|_{L^2}) \\
& &    
\|\phi_+\|_{X^{\frac{1}{2},\frac{1}{3},1}_+} + 
\|\phi_-\|_{X^{\frac{1}{2},\frac{1}{3},1}_-} \le c 
 T^{\frac{1}{6}} ( 
\|\phi_+(0)\|_{H^{\frac{1}{2}}} + \|\phi_-(0)\|_{H^{\frac{1}{2}}}) \}
\end{eqnarray*}
Taking an element $(\psi_{\pm},\phi_{\pm}) \in M_T$, the nonlinear term on the 
right hand side of (\ref{1.3}) is estimated in the 
$X^{0,\frac{1}{3},1}_{\pm}(T)$ - norm by use of Propositions \ref{linest} and 
\ref{Prop. 2.3} (we omit $T$ here and in the following)
\begin{eqnarray*}
& & \| \int_0^t e^{\mp i(t-s)|D|}
 \Pi_{\pm}(D)
  (\frac{1}{2}(\phi_+(s)+\phi_-(s)) 
  \beta(\Pi_+(D)\psi_+(s) \\ 
& & +\Pi_-(D) \psi_-(s)) ds 
   \|_{X^{0,\frac{1}{3},1}_{\pm}} \\
  & \le & c T^{\frac{1}{6}} \| \int_0^t e^{\mp i(t-s)|D|}
 \Pi_{\pm}(D)
  (\frac{1}{2}(\phi_+(s)+\phi_-(s)) 
  \beta(\Pi_+(D)\psi_+(s) + \\
  & & \Pi_-(D) \psi_-(s)) ds 
   \|_{X^{0,\frac{1}{2},1}_{\pm}} \\
   & \le & c T^{\frac{1}{3}} \| 
  \Pi_{\pm}(D)(\phi_+ +\phi_-) 
  \beta(\Pi_+(D)\psi_+  
 + \Pi_-(D) \psi_-)  
   \|_{X^{0,-\frac{1}{3},\infty}_{\pm}}\\
   & \le & c T^{\frac{1}{3}} (\| \phi_+\|_{X^{\frac{1}{2},\frac{1}{3},1}_+} + 
\| 
\phi_-\|_{X^{\frac{1}{2},\frac{1}{3},1}_-}) (\| 
\psi_+\|_{X^{0,\frac{1}{3},1}_+} 
+ \| \psi_-\|_{X^{0,\frac{1}{3},1}_+})\\
& \le & c T^{\frac{1}{2}} T^{\frac{1}{6}} (\|\phi_+(0)\|_{H^{\frac{1}{2}}} + 
\|\phi_-(0)\|_{H^{\frac{1}{2}}}) ( \|\psi_+(0)\|_{L^2} + \|\psi_-(0)\|_{L^2}) 
\\
& \le & c T^{\frac{1}{6}} ( \|\psi_+(0)\|_{L^2} + \|\psi_-(0)\|_{L^2}) \, ,
 \end{eqnarray*}
 where in the last line we used (\ref{1.6}). 
 
 The linear terms on the right hand side of (\ref{1.3}) are estimated as 
follows:
 $$ \| e^{\mp it|D|} \psi_{\pm}(0)\|_{X^{0,\frac{1}{3},1}_{\pm}} \le c 
T^{\frac{1}{6}} \| e^{\mp it|D|} \psi_{\pm}(0)\|_{X^{0,\frac{1}{2},1}_{\pm}} 
\le 
c T^{\frac{1}{6}} \|\psi_{\pm}(0)\|_{L^2} \, , $$
 and by Prop. \ref{linest} and Prop. \ref{Prop. 2.2'}:     
\begin{eqnarray*}
&& \| \int_0^t e^{\mp i(t-s)|D|} \beta \psi_{\mp}(s) ds 
\|_{X^{0,\frac{1}{3},1}_{\pm}}  \le  c T^{\frac{1}{6}} \| \int_0^t e^{\mp 
i(t-s)|D|} \beta \psi_{\mp}(s) ds \|_{X^{0,\frac{1}{2},1}_{\pm}} \\ 
&& \le c T^{\frac{1}{3}} \|\psi_{\mp}\|_{X^{0,-\frac{1}{3},\infty}_{\pm}} 
  \le  c T^{\frac{2}{3}} \|\psi_{\mp}\|_{L^2([0,T],L^2)} \\ 
&& \le  c T \|\psi_{\mp}\|_{X^{0,\frac{1}{3},1}_{\mp}} 
 \le  c T^{\frac{7}{6}} (\|\psi_+(0)\|_{L^2} + \|\psi_-(0)\|_{L^2})  \, . 
\end{eqnarray*}
Next we consider the right hand side of (\ref{1.4}).
\begin{eqnarray*}
 & & \| \int_0^t e^{\mp i(t-s)A^{\frac{1}{2}}} A^{-{\frac{1}{2}}} \langle 
\beta(\Pi_+(D) \psi_+(s)+\Pi_-(D)\psi_-(s)), \\
 & & \hspace{1cm} \Pi_+(D) \psi_+(s) + \Pi_-(D)\psi_-(s) \rangle ds 
\|_{X^{\frac{1}{2},\frac{1}{3},1}} \\
& & \le c T^{\frac{1}{6}} \| \int_0^t e^{\mp i(t-s)A^{\frac{1}{2}}} 
A^-{\frac{1}{2}} \langle \beta(\Pi_+(D) \psi_+(s)+\Pi_-(D)\psi_-(s)), \\
 & & \hspace{1cm} \Pi_+(D) \psi_+(s) + \Pi_-(D)\psi_-(s) \rangle ds 
\|_{X^{\frac{1}{2},\frac{1}{2},1}} \\
 & & \le c T^{\frac{1}{3}} \|\langle \beta(\Pi_+(D) \psi_+ +\Pi_-(D)\psi_-), 
 \Pi_+(D) \psi_+ + \Pi_-(D)\psi_- \rangle 
\|_{X^{-\frac{1}{2},-\frac{1}{3},\infty}}\\
 & & \le c T^{\frac{1}{3}} (\|\psi_+\|_{X^{0,\frac{1}{3},1}_+} + 
\|\psi_-\|_{X^{0,\frac{1}{3},1}_-})^2 \\
 & & \le c T^{\frac{1}{2}} T^{\frac{1}{6}} (\|\psi_+(0)\|_{L^2} + 
\|\psi_-(0)\|_{L^2} )^2 \\
 & & \le c T^{\frac{1}{6}}  ( \|\phi_+(0)\|_{H^{\frac{1}{2}}} + 
\|\phi_-(0)\|_{H^{\frac{1}{2}}}) \, ,
 \end{eqnarray*}
where we used (\ref{1.1}) and also (\ref{1.7}) in the last line.

The linear terms on the right hand side of (\ref{1.4}) are handled as follows:
$$ \| e^{\mp itA^{\frac{1}{2}}} 
\phi_{\pm}(0)\|_{X^{\frac{1}{2},\frac{1}{3},1}_{\pm}} \le c T^{\frac{1}{6}} \| 
e^{\mp itA^{\frac{1}{2}}} \phi_{\pm}(0)\|_{X^{\frac{1}{2},\frac{1}{2},1}_{\pm}} 
\le c T^{\frac{1}{6}} \|\phi_{\pm}(0)\|_{H^{\frac{1}{2}}} $$ 
and
\begin{eqnarray*}
\lefteqn{\| \int_0^t e^{\mp i(t-s)A^{\frac{1}{2}}} A^{-\frac{1}{2}} 
(\phi_+(s)+\phi_-(s)) ds \|_{X^{\frac{1}{2},\frac{1}{3},1}_{\pm}}} \\  
& \le & c T^{\frac{1}{6}} \| \int_0^t e^{\mp i(t-s)A^{\frac{1}{2}}} 
A^{-\frac{1}{2}} (\phi_+(s)+\phi_-(s)) ds 
\|_{X^{\frac{1}{2},\frac{1}{2},1}_{\pm}} \\
& \le & c T^{\frac{7}{6}} ( \|\phi_+(0)\|_{H^{\frac{1}{2}}} + 
\|\phi_-(0)\|_{H^{\frac{1}{2}}}) \, . 
\end{eqnarray*}  
Here we used the following estimate
\begin{eqnarray} \nonumber
\lefteqn{\| \int_0^t e^{\mp i(t-s)A^{\frac{1}{2}}} A^{-\frac{1}{2}} 
(\phi_+(s)+\phi_-(s)) ds \|_{X^{\frac{1}{2},\frac{1}{2},1}_{\pm}}} \\ \nonumber  
& \le & c T^{\frac{1}{6}} 
(\|\phi_+\|_{X^{-\frac{1}{2},-\frac{1}{3},\infty}_{\pm}} +  
\|\phi_-\|_{X^{-\frac{1}{2},-\frac{1}{3},\infty}_{\pm}}) \\ \nonumber
& \le & c T^{\frac{1}{2}} (\|\phi_+\|_{L^2_{xt}} + \|\phi_-\|_{L^2_{xt}}) \\ 
\nonumber
& \le & c T^{\frac{5}{6}} (\|\phi_+\|_{X^{0,\frac{1}{3},1}_+} +  
(\|\phi_-\|_{X^{0,\frac{1}{3},1}_-} \\ \label{!!!} 
& \le & c T (\|\phi_+(0)\|_{H^{\frac{1}{2}}} + \|\phi_-(0)\|_{H^{\frac{1}{2}}}) 
\, .
\end{eqnarray}
Altogether we have shown that the set $M_T$ is mapped into itself. Concerning 
the 
contraction property we get similarly for the difference of the right hand 
sides 
of (\ref{1.3}) applied to functions $(\psi_{\pm},\phi_{\pm}) \in M_T $ and 
$(\tilde{\psi}_{\pm},\tilde{\psi}_{\pm}) \in M_T $ in the 
$X^{0,\frac{1}{3},1}_{\pm}$ - norm an estimate by  
\begin{eqnarray*}
\lefteqn{c T^{\frac{1}{3}}  (\|\phi_+ 
-\tilde{\phi}_+\|_{X^{\frac{1}{2},\frac{1}{3},1}_+} 
+
\|\phi_- -\tilde{\phi}_-\|_{X^{\frac{1}{2},\frac{1}{3},1}_-})
 (\| \psi_+\|_{X^{0,\frac{1}{3},1}_+}
 + \| \psi_-\|_{X^{0,\frac{1}{3},1}_-}} \\
&& + \| 
\tilde{\psi}_+\|_{X^{0,\frac{1}{3},1}_+} + \| 
\tilde{\psi}_-\|_{X^{0,\frac{1}{3},1}_-} )
+ (\|\psi_+ -\tilde{\psi}_+\|_{X^{0,\frac{1}{3},1}_+} 
+ \|\psi_- -\tilde{\psi}_-\|_{X^{0,\frac{1}{3},1}_-})\cdot  \\ 
&& \cdot(\| \phi_+\|_{X^{\frac{1}{2},\frac{1}{3},1}_+} + \| 
\phi_-\|_{X^{\frac{1}{2},\frac{1}{3},1}_-} + \| 
\tilde{\phi}_+\|_{X^{\frac{1}{2},\frac{1}{3},1}_+}
+ \| \tilde{\phi}_-\|_{X^{\frac{1}{2},\frac{1}{3},1}_-}) \\
 & \le & c T^{\frac{1}{2}} [ (\|\phi_+ 
-\tilde{\phi}_+\|_{X^{\frac{1}{2},\frac{1}{3},1}_+} +
\|\phi_- -\tilde{\phi}_-\|_{X^{\frac{1}{2},\frac{1}{3},1}_-}) \cdot \\
&& \cdot (\|\psi_+(0)\|_{L^2} + \|\tilde{\psi}_+(0)\|_{L^2} + 
\|\psi_-(0)\|_{L^2}  + 
\|\tilde{\psi}_-(0)\|_{L^2})\\
&& +( \|\psi_+ -\tilde{\psi}_+\|_{X^{0,\frac{1}{3},1}_+} 
+ \|\psi_- -\tilde{\psi}_-\|_{X^{0,\frac{1}{3},1}_-}) \cdot \\
&&  \cdot (\| \phi_+(0)\|_{H^{\frac{1}{2}}} + \|\tilde 
\phi_+(0)\|_{H^{\frac{1}{2}}} +
 \| \phi_-(0)\|_{H^{\frac{1}{2}}} + \| \tilde{\phi}_-(0)\|_{H^{\frac{1}{2}}})] 
\\
& \le & \frac{1}{2}  (\|\phi_+ 
-\tilde{\phi}_+\|_{X^{\frac{1}{2},\frac{1}{3},1}_+} +
\|\phi_- -\tilde{\phi}_-\|_{X^{\frac{1}{2},\frac{1}{3},1}_-} \\
& & + \|\psi_+ -\tilde{\psi}_+\|_{X^{0,\frac{1}{3},1}_+} 
+ \|\psi_- -\tilde{\psi}_-\|_{X^{0,\frac{1}{3},1}_-}) \, ,
\end{eqnarray*} 
using (\ref{1.5}) and (\ref{1.6}) in the last line. The linear integral term in 
(\ref{1.3}) is treated easily, and the right hand side of (\ref{1.4}) can also 
be 
estimated similarly. Thus the contraction property is proved leading to a 
unique 
(local) solution.

We now show that our local solution belongs to $C^0_t(L^2_x)$. From our 
integral 
equation we get
\begin{eqnarray*}
& & \|\psi_{\pm}\|_{C^0_t L^2_x} \le c 
\|\psi_{\pm}\|_{X^{0,\frac{1}{2},1}_{\pm}}  \\
& \hspace{-0.5cm} \le & \hspace{-0.4cm} c ( \|\psi_{\pm}(0)\|_{L^2_x} + \| 
\int_0^t e^{\mp i(t-s)|D|} 
\Pi_{\pm}(D)(\phi \beta (\Pi_+(D) \psi_+ + \Pi_-(D)\psi_-))(s) ds 
\|_{X^{0,\frac{1}{2},1}_{\pm}}\\
& & + |M| \| \int_0^t e^{\mp i(t-s)|D|} \beta \psi_{\mp}(s) ds 
\|_{X^{0,\frac{1}{2},1}_{\pm}} )\\
& \hspace{-0.5cm} \le \hspace{-0.5cm} & c(  \|\psi_{\pm}(0)\|_{L^2_x} + 
T^{\frac{1}{6}}\| \Pi_{\pm}(D)(\phi 
\beta (\Pi_+(D) \psi_+ + \Pi_-(D)\psi_-))\|_{X^{0,-\frac{1}{3},\infty}}\\
& & + |M| 
 T^{\frac{5}{6}} (\|\psi_+\|_{X^{0,\frac{1}{3},1}_+} + 
\|\psi_-\|_{X^{0,\frac{1}{3},1}_-}))\\
 & \hspace{-0.5cm} \le \hspace{-0.5cm} & c (\|\psi_{\pm}(0)\|_{L^2_x} + 
T^{\frac{1}{6}} 
(\|\phi_+\|_{X^{\frac{1}{2},\frac{1}{3},1}_+} + \| 
\phi_-\|_{X^{\frac{1}{2},\frac{1}{3},1}_-})( \| 
\psi_+\|_{X^{0,\frac{1}{3},1}_+} 
\| \psi_-\|_{X^{0,\frac{1}{3},1}_-})\\
 & & + |M| T^{\frac{5}{6}} (\|\psi_+\|_{X^{0,\frac{1}{3},1}_+} + 
\|\psi_-\|_{X^{0,\frac{1}{3},1}_-}))
\end{eqnarray*}
Here we used the estimate
$$
\| \int_0^t e^{\mp i(t-s)|D|} \beta \psi_{\mp}(s) ds 
\|_{X^{0,\frac{1}{2},1}_{\pm}} 
 \le  c T^{\frac{1}6} (\| \psi_-\|_{X^{0,-\frac{1}{3},\infty}_{\pm}} + \| 
\psi_+\|_{X^{0,-\frac{1}{3},\infty}_{\pm}}) $$
 $$
\le c T^{\frac{1}{2}}(\|\psi_-\|_{L^2_{xt}} + \psi_+\|_{L^2_{xt}})  \le  c 
T^{\frac{5}{6}} (\| \psi_+\|_{X^{0,\frac{1}{3},1}_+} + \| 
\psi_-\|_{X^{0,\frac{1}{3},1}_-})
$$
by Prop. \ref{linest} and Prop. \ref{Prop. 2.2'}.
We have shown that $\psi_{\pm} \in C^0_t(L^2_x)$. 

Next we estimate $\|\phi_{\pm}(t)\|_{H^{\frac{1}{2}}_x}$ for $0 \le t \le T$ by 
our integral equation (\ref{1.4}).
\begin{eqnarray} \nonumber
& &\|\phi_{\pm}(t)\|_{H^{\frac{1}{2}}_x}  \\ \nonumber
& \le & \|\phi_{\pm}(0)\|_{H^{\frac{1}{2}}_x} + c \|\int_0^t e^{\mp 
i(t-s)A^{\frac{1}{2}}} A^{-\frac{1}{2}} \langle \beta(\Pi_+(D)\psi_+(s) + 
\Pi_-(D)\psi_-(s)), \\ \nonumber
& & \Pi_+(D)\psi_+(s) + \Pi_-(D)\psi_-(s))\rangle 
ds\|_{X^{\frac{1}{2},\frac{1}{2},1}_{\pm}} \\ \nonumber
& & +c|m+1| \|\int_0^t e^{\mp i(t-s)A^{\frac{1}{2}}} A^{-\frac{1}{2}} \phi(s) 
ds 
\|_{X^{\frac{1}{2},\frac{1}{2},1}_{\pm}} \\ \nonumber
& \le & \|\phi_{\pm}(0)\|_{H^{\frac{1}{2}}_x} + cT^{\frac{1}{6}} \|\langle 
\beta(\Pi_+(D)\psi_+(s) + \Pi_-(D)\psi_-(s)),  \Pi_+(D)\psi_+(s) \\ \nonumber
& & + \Pi_-(D)\psi_-(s))\rangle\|_{X^{-\frac{1}{2},-\frac{1}{3},\infty}_{\pm}} 
+ 
c|m+1| \|\int_0^t e^{\mp i(t-s)A^{\frac{1}{2}}} A^{-\frac{1}{2}} \phi(s) ds 
\|_{X^{\frac{1}{2},\frac{1}{2},1}_{\pm}}\\ \nonumber
& \le & \|\phi_{\pm}(0)\|_{H^{\frac{1}{2}}_x} + c T^{\frac{1}{6}} 
(\|\psi_+\|_{X^{0,\frac{1}{3},1}_+}^2 + \|\psi_-\|_{X^{0,\frac{1}{3},1}_-}^2)\\ 
\nonumber
& & + c|m+1| \|\int_0^t e^{\mp i(t-s)A^{\frac{1}{2}}} A^{-\frac{1}{2}} \phi(s) 
ds \|_{X^{\frac{1}{2},\frac{1}{2},1}_{\pm}}\\ \nonumber
& \le & \|\phi_{\pm}(0)\|_{H^{\frac{1}{2}}_x} + cT^{\frac{1}{2}} 
(\|\psi_+(0)\|_{L^2}^2 + \|\psi_-(0)\|_{L^2}^2)\\ \label{****}
& &  + c|m+1| T (\|\phi_+(0)\|_{H^{\frac{1}{2}}} + 
\|\phi_-(0)\|_{H^{\frac{1}{2}}}) \, . 
\end{eqnarray}
Here we used (\ref{!!!}).
By our choice (\ref{1.6}) of $T$ we arrive at (\ref{*****}). 
This proves that $ \phi_{\pm} \in C^0([0,T],H^{\frac{1}{2}})$.
For our global result it is important to notice that no implicit constant 
appears in front of the first term on the right hand side.

The additional claim of the proposition  is easily proven by the contraction 
mapping principle in a similar way so that the proof of Prop. \ref{Prop. 2.4} 
is complete. 

The proof of Theorem \ref{Theorem} can now be given along the lines of the 
paper 
of Colliander-Holmer-Tzirakis \cite{CHT} for the 1D Zakharov system.\\[0.2cm]
{\bf Proof} of Theorem \ref{Theorem}.
We start by using the addition in Prop. \ref{Prop. 2.4} leading to a local 
solution with the required regularity 
properties. Because $\|\psi(t)\|_{L^2}$ is conserved we get by iteration a 
global solution if also $\|\phi(t)\|_{H^{\frac{1}{2}}}$ remains bounded. 
Otherwise 
we use our Prop. \ref{Prop. 2.4} and remark first that 
$$ \|\psi(t)\|_{L^2}^2 = \|\psi_+(t)\|_{L^2}^2 + \|\psi_-(t)\|_{L^2}^2 $$
is still conserved. This conservation law can be applied because $\psi_{\pm} 
\in 
C^0([0,T],L^2_x) \hspace{-0.6pt}.$ Without loss of generality we can now 
suppose that at some 
time t we have
$$ \|\phi_+(t)\|_{H^{\frac{1}{2}}} + \|\phi_-(t)\|_{H^{\frac{1}{2}}} >> 
\|\psi_+(t)\|_{L^2}^2 + \|\psi_-(t)\|_{L^2}^2 \, . $$
Take this time t as initial time $t=0$ so that
\begin{equation}
\label{******}
\|\phi_+(0)\|_{H^{\frac{1}{2}}} + \|\phi_-(0)\|_{H^{\frac{1}{2}}} >> 
\|\psi_+(0)\|_{L^2}^2 + \|\psi_-(0)\|_{L^2}^2 \, .
\end{equation}
Then (\ref{1.7}) is automatically satisfied. We define
$$ T \sim \frac{1}{(\|\phi_+(0)\|_{H^{\frac{1}{2}}} + 
\|\phi_-(0)\|_{H^{\frac{1}{2}}})^2} \, , $$  
so that (\ref{1.5}) and (\ref{1.6}) are fulfilled.
From our estimate (\ref{*****}) we conclude that it is possible to use the 
local 
existence result $l$ times with time intervals of length $T,$ before the 
quantity $ \|\phi_+(t)\|_{H^{\frac{1}{2}}} + \|\phi_-(t)\|_{H^{\frac{1}{2}}} $ 
doubles. Here we have
$$ l \sim \frac{\|\phi_+(0)\|_{H^{\frac{1}{2}}} + 
\|\phi_-(0)\|_{H^{\frac{1}{2}}}}{T^{\frac{1}{2}}(\|\psi_+(0)\|_{L^2}^2 + 
\|\psi_-(0)\|_{L^2}^2 +1)} \, . $$
  
After these $l$ iterations we arrive at the time
$$ lT \sim \frac{\|\phi_+(0)\|_{H^{\frac{1}{2}}} + 
\|\phi_-(0)\|_{H^{\frac{1}{2}}}}{\|\psi_+(0)\|_{L^2}^2 + \|\psi_-(0)\|_{L^2}^2 
+1} T^{\frac{1}{2}} \sim  \frac{1}{\|\psi_+(0)\|_{L^2}^2 + 
\|\psi_-(0)\|_{L^2}^2 
+1} \, . $$
This quantity is independent of $ \|\phi_+(0)\|_{H^{\frac{1}{2}}} + 
\|\phi_-(0)\|_{H^{\frac{1}{2}}} $. Using conservation of $\|\psi_+(t)\|_{L^2}^2 
+ \|\psi_-(t)\|_{L^2}^2$ it is thus possible to repeat the whole procedere with 
time steps of equal length. This proves the global existence result.\\[0.2cm]
{\bf Proof} of Theorem \ref{Theorem'}. By the Leibniz rule for fractional 
derivatives from (\ref{1.1}),(\ref{1.2}) one easily gets the following 
estimates for the nonlinearities for arbitrary $s \ge 0$ :
\begin{eqnarray*}
\lefteqn{ \| \langle \beta \Pi_{\pm 1}(D)\psi, \Pi_{\pm 2}(D)\psi ' \rangle 
\|_{X^{s-\frac{1}{2},-\frac{1}{3},\infty}_{\pm 3}}} \\
& \le & c(\|\psi\|_{X^{0,\frac{1}{3},1}_{\pm 1}} \|\psi 
'\|_{X^{s,\frac{1}{3},1}_{\pm 2}} + \|\psi\|_{X^{s,\frac{1}{3},1}_{\pm 1}} 
\|\psi '\|_{X^{0,\frac{1}{3},1}_{\pm 2}})
\end{eqnarray*}
and
\begin{eqnarray*}
\lefteqn{ \| \Pi_{\pm 2}(D)(\phi \beta \Pi_{\pm 1}(D)\psi) 
\|_{X^{s,-\frac{1}{3},\infty}_{\pm 2}}} \\
& \le & c(\|\phi\|_{X^{s+\frac{1}{2},\frac{1}{3},1}_{\pm 3}} 
\|\psi\|_{X^{0,\frac{1}{3},1}_{\pm 1}} + 
\|\phi\|_{X^{\frac{1}{2},\frac{1}{3},1}_{\pm 3}} 
\|\psi\|_{X^{s,\frac{1}{3},1}_{\pm 1}})\, .
\end{eqnarray*}
These estimates allow to construct a local solution 
with the required properties by the contraction mapping principle with  an 
existence time $$ T \sim \frac{1}{\|\psi_+(0)\|_{L^2} + \|\psi_-(0)\|_{L^2} + 
\|\phi_+(0)\|_{H^{\frac{1}{2}}} + \|\phi_-(0)\|_{H^{\frac{1}{2}}}}   $$
similarly as in the proof of Prop. \ref{Prop. 2.4}. By uniqueness the global 
solution of Theorem \ref{Theorem} coincides locally with this solution. Thus 
the claim of Theorem \ref{Theorem'} follows. 
\section{The estimates for the nonlinearities}
In this section we give the proof of Prop. \ref{Prop. 2.3}. We first show that 
the estimates (\ref{1.1}) and (\ref{1.2}) are completely equivalent to each 
other. By duality (\ref{1.2}) is equivalent to the estimate
\begin{equation}
\label{2.1}
\left| \int \langle \Pi_{\pm 2}(D) (\phi \beta \Pi_{\pm 1}(D)\psi),\psi' 
\rangle 
dx dt \right| \le c \|\phi\|_{X^{\frac{1}{2},\frac{1}{3},1}_{\pm 3}} 
\|\psi\|_{X^{0,\frac{1}{3},1}_{\pm 1}} \|\psi'\|_{X^{0,\frac{1}{3},1}_{\pm 2}} 
\, .
 \end{equation}
The left hand side equals
\begin{equation}
\label{2.1'} \left| \int \phi \langle \beta \Pi_{\pm 1}(D)\psi,\Pi_{\pm 2}\psi' 
\rangle dx dt \right| \, . 
\end{equation}
Thus (\ref{2.1}) is equivalent to (\ref{1.1}).

The complete null structure of the system detected by d'Ancona, Foschi and 
Selberg has the following consequences (cf. \cite{AFS}). 
Denoting
$$ \sigma_{\pm 1,\pm 2}(\eta,\zeta) := \Pi_{\pm 2}(\zeta) \beta \Pi_{\pm 
1}(\eta) = \beta \Pi_{\mp 2}(\zeta) \Pi_{\pm 1}(\eta) \, , $$
we remark that by orthogonality this quantity vanishes if $\pm_1 \eta$ and 
$\pm_2  \zeta$ line up in the same direction whereas in general (cf. 
\cite{AFS1}, Lemma 1):
\begin{lemma}
\label{Lemma 3.1}
$$\sigma_{\pm 1,\pm 2}(\eta,\zeta) = O(\angle (\pm_1 \eta, \pm_2 \zeta)) \, ,$$ 
where $\angle(\eta,\zeta)$ denotes the angle between the vectors $\eta$ and 
$\zeta$.
\end{lemma}
Consequently we get
\begin{eqnarray} \nonumber
\lefteqn{ |\langle \beta \Pi_{\pm 1}(D)\psi ,\Pi_{\pm 2}(D)\psi ' 
  \rangle ^{\tilde{}}(\tau,\xi)|} \\ \nonumber
& \le & \int |\langle \beta \Pi_{\pm 1}(\eta)\tilde \psi(\lambda,\eta), 
\Pi_{\pm 
2}(\eta - \xi) \tilde{\psi}'(\lambda - \tau,\eta - \xi)\rangle| d\lambda d\eta 
\\ \label{2.2}
& = & 
\int |\langle \Pi_{\pm 2}(\eta - \xi) \beta \Pi_{\pm 1}(\eta)\tilde 
\psi(\lambda,\eta),  \tilde{\psi}'(\lambda - \tau,\eta - \xi)\rangle| d\lambda 
d\eta\\ \nonumber
& \le & c \int \Theta_{\pm 1,\pm 2} \, |\tilde{\psi}(\lambda,\eta)|\, 
|\tilde{\psi}'(\lambda - \tau, \eta - \xi)| d\lambda d\eta \, ,
\end{eqnarray} 
where $\Theta_{\pm 1,\pm 2} = \angle(\pm_1 \eta,\pm_2(\eta - \xi))$.

We also need the following elementary estimates which can be found in 
\cite{AFS}, 
section 5.1.
\begin{lemma}
\label{Lemma 3.2}
Denoting
$$ A_{\pm 1} = \tau \pm_1 |\xi| \, , \, B = \lambda + |\eta| \, , \, C_{\pm} = 
\lambda - \tau \pm |\eta - \xi| \, , \, \Theta_{\pm} = 
\angle(\eta,\pm(\eta-\xi)) $$
and
$$ \rho_+ = |\xi| - \left| |\eta| - |\eta-\xi| \right| \, , \, \rho_- = 
|\eta|+|\eta-\xi| - |\xi | $$
the following estimates hold:
$$ \Theta_+^2 \sim \frac{|\xi| \rho_+}{|\eta| |\eta - \xi|} \, , \, \Theta_-^2 
\sim \frac{(|\eta|+|\eta-\xi|)\rho_-}{|\eta| |\eta-\xi|} \sim 
\frac{\rho_-}{\min(|\eta|,|\eta - \xi|)} $$ 
as well as
$$ \rho_{\pm} \le 2 \min(|\eta|,|\eta-\xi|) $$
and 
$$ \rho_ {\pm} \le |A_{\pm 1}| + |B| + |C_{\pm}| \, . $$
\end{lemma} 
{\bf Proof:}
We only prove the last estimate. We have
\begin{eqnarray*}
\rho_+ & \le & |\xi| \pm |\eta| \mp|\eta-\xi| = |\xi| \mp \tau \pm \lambda \pm 
|\eta| \pm \tau \mp \lambda \mp|\eta - \xi| \\
& \le & \left| |\xi| \mp \tau\right| + \left|\lambda + |\eta| \right| + \left| 
\tau - \lambda - |\eta-\xi| \right| \\
& \le & |A_{\mp}| + |B| + |C_+|
\end{eqnarray*}
and
\begin{eqnarray*}
\rho_- & = & (\lambda + |\eta|)+(\tau-\lambda+|\eta-\xi|) - (\tau + |\xi|) \\
& \le & |\lambda + |\eta|| + |\lambda -\tau-|\eta-\xi|| + |\tau + |\xi|| \\
& = & |B| + |C_-| + |A_+|
\end{eqnarray*}
as well as for $\tau \ge 0$:
$$ \rho_- \le \lambda + |\eta| + \tau-\lambda + |\eta-\xi| \le |B| + |C_-| $$
and for $\tau \le 0$:
$$ \rho_- \le |\lambda+|\eta||+|\tau-\lambda+|\eta-\xi|| + |\tau| + |\xi| \le 
|B| + |C_-| + |A_-| \, . $$ 
{\bf Proof of Prop. \ref{Prop. 2.3}:}
In order to prove (\ref{2.1}) first for the signs $\pm_1 = +$ and $\pm_2 = \pm$ 
and taking into account (\ref{2.1'}) and (\ref{2.2}) we have to show (recalling 
$\Theta_{\pm} := \angle(\eta,\pm(\eta-\xi))$):  
\begin{eqnarray}
\nonumber
I_{\pm} & := & \left| \int \int \Theta_{\pm} \tilde{\psi}(\lambda,\eta) 
\tilde{\psi}'(\lambda-\tau,\eta-\xi) d\lambda d\eta \tilde{\phi}(\tau,\xi) 
d\tau 
d\xi \right| \\
\label{2.3}
& \le & c \|\psi\|_{X^{0,\frac{1}{3},1}_+} \|\psi'\|_{X^{0,\frac{1}{3},1}_+} 
\|\phi\|_{X^{\frac{1}{2},\frac{1}{3},1}_{\pm 1}} \, .
\end{eqnarray}  
We may assume here without loss of generality that the Fourier transforms are 
nonnegative. Defining
\begin{eqnarray*}
\tilde{F}(\lambda,\eta) &:= &\langle \lambda + |\eta| \rangle^{\frac{1}{3}} 
\tilde{\psi}(\lambda,\eta) \\
\tilde{G_{\pm}}(\lambda,\eta) &:= &\langle \lambda \pm |\eta| 
\rangle^{\frac{1}{3}} \tilde{\psi}'(\lambda,\eta) \\
\tilde{H}_{\pm}(\tau,\xi) &:= &\langle \tau \pm |\xi| \rangle^{\frac{1}{3}} 
\langle \xi \rangle^{\frac{1}{2}} \tilde{\phi}(\tau,\xi) 
\end{eqnarray*}
we thus have to show
\begin{eqnarray*}
J_{\pm} & := & \left| \int \int  \Theta_{\pm} 
\frac{\tilde{F}(\lambda,\eta)}{\langle B \rangle ^{\frac{1}{3}}}
\frac{\tilde{G}_{\pm}(\lambda-\tau,\eta-\xi)}{\langle C_{\pm} \rangle 
^{\frac{1}{3}}}
\frac{\tilde{H}_{\pm 1}(\tau,\xi)}{\langle A_{\pm 1} \rangle ^{\frac{1}{3}} 
\langle \xi \rangle^{\frac{1}{2}}} d\lambda d\eta d\tau d\xi \right| \\
& \le & c \|F\|_{X^{0,0,1}_+} \|G_{\pm}\|_{X^{0,0,1}_{\pm}} \|H_{\pm 
1}\|_{X^{0,0,1}_{\pm 1}} \, . 
\end{eqnarray*}  
 Let us first consider the low-frequency case, where $\min(|\eta|,|\eta-\xi|) 
\le 1$. Assuming without loss of generality (by symmetry) $|\eta| \le 1$ we 
estimate
\begin{eqnarray*}
I_{\pm} & \le & \|\psi\|_{L_t^4 (L_x^{\infty})} \|\psi'\|_{L^4_t(L^2_x)} 
\|\phi\|_{L^2_tL^2_x} \\
& \le & \|\psi\|_{X_+^{2,\frac{1}{4}+}} \|\psi'\|_{X_{\pm}^{0,\frac{1}{4}+}} 
\|\phi\|_{L^2_tL^2_x} \\
& \le & \|\psi\|_{X_+^{0,\frac{1}{4}+}} \|\psi'\|_{X_{\pm}^{0,\frac{1}{4}+}}  
\|\phi\|_{X^{0,0}} \, ,
\end{eqnarray*}   
which implies the desired estimate. From now on we assume $|\eta|,|\eta-\xi| 
\ge 
1$. \\
{\bf Estimate for $J_+$:} We use
$$ \Theta_+ \le \frac{|\xi|^{\frac{1}{2}} 
\rho_+^{\frac{1}{2}}}{|\eta|^{\frac{1}{2}}  |\eta - 
\xi|^{\frac{1}{2}}} \le c \frac{\langle \xi \rangle^{\frac{1}{2}} 
\rho_+^{\frac{1}{6}}}{\langle \eta \rangle^{\frac{1}{2}} \langle \eta - \xi 
\rangle^{\frac{1}{2}}} (\langle A_ {\pm} \rangle^{\frac{1}{3}} +  \langle B  
\rangle^{\frac{1}{3}} + \langle C_+ \rangle^{\frac{1}{3}})
$$
and also
$$ \rho_+ \le 2 \min(|\eta|,|\eta - \xi|) 
\, . $$
We thus get
$$ J_+ \le c(I_1^+ + I_2^+ + I_3^+) \, , $$
where  
$$
I_1^+  =  \int \frac{\tilde{F}(\lambda,\eta)}{\langle \eta 
\rangle^{\frac{5}{12}} \langle B \rangle^{\frac{1}{3}}} 
\frac{\tilde{G}_+(\lambda-\tau,\eta-\xi)}{\langle \eta - \xi 
\rangle^{\frac{5}{12}} \langle C_+ \rangle^{\frac{1}{3}}} \tilde{H}_{\pm 
1}(\tau,\xi) d\lambda d\eta d\tau d\xi \, , $$ 
$$
I_2^+  =  \int \frac{\tilde{F}(\lambda,\eta)}{\langle \eta 
\rangle^{\frac{1}{3}}} \frac{\tilde{G}_+(\lambda-\tau,\eta-\xi)}{\langle \eta - 
\xi \rangle^{\frac{1}{2}} \langle C_+ \rangle^{\frac{1}{3}}} 
\frac{\tilde{H}_{\pm 1}(\tau,\xi)}{\langle A \rangle^{\frac{1}{3}}} d\lambda 
d\eta d\tau d\xi \, , $$
$$
I_3^+  =  \int \frac{\tilde{F}(\lambda,\eta)}{\langle \eta 
\rangle^{\frac{1}{2}}\langle B \rangle^{\frac{1}{3}}} 
\frac{\tilde{G}_+(\lambda-\tau,\eta-\xi)}{\langle \eta - \xi 
\rangle^{\frac{1}{3}}} \frac{\tilde{H}_{\pm 1}(\tau,\xi)}{\langle A 
\rangle^{\frac{1}{3}}} d\lambda d\eta d\tau d\xi \, . $$ 
We only consider $I_1^+$ and $I_2^+$, because $I_3^+$  is similar to $I_2^+$.\\
{\bf Estimate for $I_1^+$ :}
H\"older's inequality and Parseval's identity give
$$
I_1^+  \le  c \|H_{\pm 1}\|_{L^2_{xt}} \left\| {\cal F}^{-1}( 
\frac{\tilde{F}(\lambda , \eta)}{\langle \eta \rangle^{\frac{5}{12}} \langle B 
\rangle^{\frac{1}{3}}})\right\|_{L^4_{xt}} \left\| {\cal F}^{-1}( 
\frac{\tilde{G_+}(\lambda -\tau, \eta - \xi)}{\langle \eta - \xi 
\rangle^{\frac{5}{12}} \langle C_+ \rangle^{\frac{1}{3}}})\right\|_{L^4_{xt}} 
\, 
. $$
Concerning the last two factors we use Strichartz' inequality for the wave 
equation which gives for $U(t) = e^{it|D|}$:
$$ \|U(t)u_0\|_{L^8_t (H^{-\frac{5}{8},8}_x)} \le c \|u_0\|_{L^2} \, . $$ 
This implies by Prop. \ref{Prop. 1.1}:
$$ \|f\|_{L^8_t (H^{-\frac{5}{8},8}_x)} \le c 
\|U(-t)f\|_{B^{\frac{1}{2}}_{2,1}L^2_x} 
 = c \|f\|_{X^{0,\frac{1}{2},1}_+}
$$ 
Moreover we have 
$$ \|f\|_{L^2_t L^2_x} = \|U(-t)f\|_{L^2_t L^2_x} \le c \|U(-t)f\|_{B^0_{2,1} 
L^2_x} = c \|f\|_{X^{0,0,1}_+} \, . $$
Complex interpolation gives by \cite{BL}, Thm. 6.4.5:
$$\|f\|_{L^4_t (H^{-\frac{5}{12},4}_x)} \le c 
\|U(-t)f\|_{B^{\frac{1}{3}}_{2,1}L^2_x} 
 = c \|f\|_{X^{0,\frac{1}{3},1}_+}
$$ 
This is equivalent to 
$$\|f\|_{L^4_t L^4_x} \le c \|U(-t)f\|_{B^{\frac{1}{3}}_{2,1} 
H^{\frac{5}{12}}_x} 
 = c \|f\|_{X^{\frac{5}{12},\frac{1}{3},1}_+} \, .
$$       
Thus we get
$$ I_1^+ \le c \|H_{\pm 1}\|_{X^{0,0,1}_{\pm 1}} \|F\|_{X^{0,0,1}_+} 
\|G_+\|_{X^{0,0,1}_+} \, , $$
where we used the embedding $ X^{0,0,1}_{\pm} \subset L^2_{xt}$.\\
{\bf Estimate for $I_2^+$:}
Using Parseval's identity and H\"older's inequality we get
$$ I_2^+ \hspace{-0.1cm} \le \hspace{-0.1cm} c \left\| {\cal 
F}^{-1}(\frac{F(\lambda,\eta)}{\langle \eta 
\rangle^{\frac{1}{3}}})\right\|_{L^2_t( L^3_x)}  \left\| {\cal 
F}^{-1}(\frac{G_+(\lambda-\tau,\eta- \xi)}{\langle \eta - \xi 
\rangle^{\frac{1}{2}} \langle C_+ \rangle^{\frac{1}{3}}})\right\|_{L^3_t( 
L^6_x)}
 \left\| {\cal F}^{-1}(\frac{H_{\pm 1}}{\langle A 
\rangle^{\frac{1}{3}}})\right\|_{L^6_t (L^2_x)} \hspace{-0.1cm}. $$
The first factor is estimated using Sobolev's embedding theorem by 
$\|F\|_{L^2_{xt}}$.\\
Concerning the last factor we estimate by Sobolev and Minkowski's inequality as 
follows:
\begin{eqnarray}
\label{!}
\|f\|_{L^6_t(L^2_x)} & = & \|U(\mp t)f\|_{L^6_t(L^2_x)} \le \|U(\mp 
t)f\|_{L^2_x(L^6_t)} \\ \nonumber
\le c \|U(\mp t)f\|_{L^2_x H^{\frac{1}{3}}_t} & = & \|U(\mp 
t)f\|_{H^{\frac{1}{3}}_t L^2_x} = \|f\|_{X^{0,\frac{1}{3}}_{\pm}} \, .
\end{eqnarray}
Thus the last factor can be estimated by $ c\|H_{\pm 1}\|_{L^2_{xt}} \le 
c\|H_{\pm 
1}\|_{X^{0,0,1}_{\pm 1}}$.\\
Concerning the second factor we start with Strichartz' estimate
$$ \|U(t)u_0\|_{L^4_t (B^{-\frac{3}{4}}_{\infty,2})} \le c \|u_0\|_{L^2_x}\, , 
$$
which implies by Prop. \ref{Prop. 1.1}
$$ \|f\|_{L^4_t (B^{-\frac{3}{4}}_{\infty,2})} \le c 
\|U(-t)f\|_{B^{\frac{1}{2}}_{2,1} L_x^2} = c \|f\|_{X^{0,\frac{1}{2},1}_+} \, . 
$$
Moreover we have
$$ \|f\|_{L^2_t (B^0_{2,2})} = \|f\|_{X^{0,0}} = \|U(-t)f\|_{L^2_tL^2_x} \le c 
\|U(-t)f\|_{B^0_{2,1} L^2_x} = \|f\|_{^{0,0,1}_+} \, . $$
We now use the complex interpolation method. By \cite{BL}, Thm. 6.4.5 we have
$$ (B^{-\frac{3}{4}}_{\infty,2},B^0_{2,2})_{[\frac{2}{3}]} = 
B^{-\frac{1}{2}}_{6,2} \quad {\mbox {and also}} \quad 
(X_+^{0,\frac{1}{2},1},X_+^{0,0,1})_{[\frac{2}{3}]} = X_+^{0,\frac{1}{3},1} \, 
,$$
so that we get with $B^{-\frac{1}{2}}_{6,2} \subset H^{-\frac{1}{2},6}$ 
(\cite{BL}, Thm. 6.4.4)
$$ \|f\|_{L_t^3 (H^{-\frac{1}{2},6})} \le c \|f\|_{L^3_t 
(B^{-\frac{1}{2}}_{6,2})} \le c \|U(-t)f\|_{B^{\frac{1}{3}}_{2,1} L^2_x} = 
\|f\|_{X^{0,\frac{1}{3},1}_+} \, , $$
which implies
\begin{equation}
\label{!!}
\|f\|_{L^3_t (L^6_x)} \le c \|f\|_{X^{\frac{1}{2},\frac{1}{3},1}_+} \, . 
\end{equation}
Thus the second factor is estimated by $\|G_+\|_{X^{0,0,1}_+}$.\\
{\bf Estimate for $J_-$:}
If $ |\eta| << |\eta - \xi|$ we have $|\xi| \sim |\eta - \xi|$, thus by Lemma 
\ref{Lemma 3.2}:
$$\Theta_-^2 \sim \frac{\rho_-}{\min(|\eta|,|\eta - \xi|)} \sim \frac{|\xi| 
\rho_-}{|\eta| |\eta-\xi|} \, , $$
so that
$$ \Theta_- \le c \frac{\langle \xi \rangle ^{\frac{1}{2}} 
\rho_-^{\frac{1}{6}}}{\langle \eta \rangle^{\frac{1}{2}} \langle \eta - \xi 
\rangle^{\frac{1}{2}}} (\langle A_{\pm} \rangle^{\frac{1}{3}} + \langle B 
\rangle^{\frac{1}{3}} + \langle C_- \rangle^{\frac{1}{3}}) \, . $$
Because also $ \rho_- \le 2 \min(|\eta|,|\eta - \xi|)$ the same estimates as 
for 
$J_+$ can be given. If $|\eta| >> |\eta-\xi|$, we have $|\xi| \ge 
||\eta|-|\eta-\xi|| \sim |\eta|$ and the same estimate for $\Theta_-$ holds. 
This is also true if $|\xi| \sim |\eta| \sim |\eta-\xi|$.

It remains to consider $J_-$ in the case $|\xi| << |\eta| \sim |\eta-\xi|$, 
which we assume from now on. We then have
$$ \Theta_- \le  \frac{\rho_-^{\frac{1}{2}}}{\langle \eta \rangle^{\frac{1}{4}} 
\langle \eta - \xi \rangle^{\frac{1}{4}}} $$
and thus
\begin{eqnarray*}
J_- & \le & c \left| \int \int \rho_-^{\frac{1}{2}} 
\frac{\tilde{F}(\lambda,\eta)}{\langle \eta \rangle^{\frac{1}{4}} \langle B 
\rangle^{\frac{1}{3}}} \frac{\tilde{G}_-(\lambda - \tau, \eta - \xi)}{\langle 
\eta - \xi \rangle^{\frac{1}{4}} \langle C_-\rangle^{\frac{1}{3}}} 
\frac{\tilde{H}_{\pm 1}(\tau,\xi)}{\langle \xi \rangle ^{\frac{1}{2}} \langle 
A_ 
{\pm 1} \rangle^{\frac{1}{3}}} d\lambda d\eta d\tau d\xi \right| \, .
\end{eqnarray*}   
Using the estimates $\rho_- \le 2 \min(|\eta|,|\eta-\xi|)$ and $ \rho_- \le 
|A_{\pm 1}| +|B| + |C_-|$ (cf. Lemma \ref{Lemma 3.2}) we get  
$$ J_- \le c(I_1^- + I_2^- + I_3^-) \, , $$
where
\begin{eqnarray*}
I_1^- & = & \int \int \frac{\tilde{F}(\lambda,\eta)}{\langle \eta 
\rangle^{\frac{1}{6}} \langle B \rangle^{\frac{1}{3}}} 
\frac{\tilde{G}_-(\lambda 
- \tau, \eta - \xi)}{\langle \eta - \xi \rangle^{\frac{1}{6}} \langle 
C_-\rangle^{\frac{1}{3}}} \frac{\tilde{H}_{\pm 1}(\tau,\xi)}{\langle \xi 
\rangle 
^{\frac{1}{2}}} d\lambda d\eta d\tau d\xi \\
I_2^- & = & \int \int \frac{\tilde{F}(\lambda,\eta)}{\langle \eta 
\rangle^{\frac{1}{12}}} \frac{\tilde{G}_-(\lambda - \tau, \eta - \xi)}{\langle 
\eta - \xi \rangle^{\frac{1}{4}} \langle C_-\rangle^{\frac{1}{3}}} 
\frac{\tilde{H}_{\pm 1}(\tau,\xi)}{\langle \xi \rangle ^{\frac{1}{2}} \langle 
A_{\pm 1} \rangle^{\frac{1}{3}}} d\lambda d\eta d\tau d\xi \\ 
I_3^- & = & \int \int \frac{\tilde{F}(\lambda,\eta)}{\langle \eta 
\rangle^{\frac{1}{4}} \langle B \rangle^{\frac{1}{3}}} 
\frac{\tilde{G}_-(\lambda 
- \tau, \eta - \xi)}{\langle \eta - \xi \rangle^{\frac{1}{12}}} 
\frac{\tilde{H}_{\pm 1}(\tau,\xi)}{\langle \xi \rangle ^{\frac{1}{2}} \langle 
A_{\pm 1} \rangle^{\frac{1}{3}}} d\lambda d\eta d\tau d\xi \, . 
\end{eqnarray*}  
The terms $I_2^-$ and $I_3^-$ are similar, so that we concentrate on $I_1^-$ 
and 
$I_2^-$.\\
{\bf Estimate for $I_1^-$:}
We have
\begin{eqnarray*}
I_1^{-} &
 \le & \| \tilde{H}_{\pm 1}\|_{L^2_{\xi \tau}}
\left\| \int \langle \xi \rangle^{-\frac{1}{2}} 
\frac{\tilde{F}(\lambda,\eta)}{\langle \eta \rangle^{\frac{1}{6}} \langle B 
\rangle^{\frac{1}{3}}} \frac{\tilde{G}_-(\lambda - \tau, \eta - \xi)}{\langle 
\eta - \xi \rangle^{\frac{1}{6}} \langle C_-\rangle^{\frac{1}{3}}}  d\lambda 
d\eta \right\|_{L^2_{\xi \tau}} \\
&  = & \|H_{\pm 1}\|_{L^2_{xt}}
\left\| \int \langle \xi \rangle^{-\frac{1}{2}} 
\frac{\tilde{F}(\lambda,\eta)}{\langle \eta \rangle^{\frac{1}{6}} \langle 
\lambda + |\eta| \rangle^{\frac{1}{3}}} \frac{\tilde{G}'(\tau - \lambda, \xi - 
\eta)}{\langle \xi - \eta \rangle^{\frac{1}{6}} \langle \tau - \lambda + |\xi - 
\eta|\rangle^{\frac{1}{3}}}  d\lambda d\eta \right\|_{L^2_{\xi \tau}}  , 
\end{eqnarray*}  
where $ \tilde{G}'(\lambda,\eta) := \tilde{G}_-(-\lambda,-\eta)$. This shows 
that we in fact are in the (+,+)-case. We also remark that we assumed $|\xi| << 
|\eta|  \sim  |\xi - \eta|$. Using Prop. \ref{Prop. 3.2} we arrive at
$$ I_1^- \le c \|H_{\pm 1}\|_{L^2_{xt}} \|F\|_{X^{0,0,1}_+} 
\|G'\|_{X^{0,0,1}_+} 
\le c \|H_{\pm 1}\|_{X^{0,0,1}_{\pm 1}} \|F\|_{X^{0,0,1}_+} 
\|G_-\|_{X^{0,0,1}_-} \, . $$ 
{\bf Estimate for $I_2^-$:}
Parseval's identity and H\"older's inequality imply
\begin{eqnarray*}
I_2^- & \le & c
\|{\cal F}^{-1}(\frac{\tilde{F}}{\langle \eta \rangle^{\frac{1}{12}}})\|_{L^2_t 
(L^{\frac{24}{11}}_x)} 
\|{\cal F}^{-1}(\frac{\tilde{G}_-}{\langle \eta - \xi \rangle^{\frac{1}{4}} 
\langle C_-\rangle^{\frac{1}{3}}})\|_{L^6_t (L^{\frac{8}{3}}_x)} \\
& & \|{\cal F}^{-1}(\frac{\tilde{H}_{\pm 1}}{\langle \xi \rangle^{\frac{1}{2}} 
\langle A_{\pm 1}\rangle^{\frac{1}{3}}})\|_{L^3_t (L^6_x)} \, . 
\end{eqnarray*}  
The first factor is controlled using Sobolev's embedding $H^{\frac{1}{12}} 
\subset L^{\frac{24}{11}}$ by $\|F\|_{L^2_{xt}}$, the last factor is handled as 
before using the estimate (\ref{!!}), and the second one similarly as before as 
follows. First, Sobolev's embedding in $x$ gives
$$ \|{\cal F}^{-1}(\frac{\tilde{G}_-}{\langle \eta - \xi \rangle^{\frac{1}{4}} 
\langle C_-\rangle^{\frac{1}{3}}})\|_{L^6_t (L^{\frac{8}{3}}_x)} \le c \|{\cal 
F}_t^{-1} (\frac{\widehat{G}_-}{\langle C_- 
\rangle^{\frac{1}{3}}})\|_{L^6_t(L_x^2)} \, . $$ 
Now we use (\ref{!}) so that the second factor is estimated by 
$\|G_-\|_{L^2_{xt}} \le c \|G_-\|_{X^{0,0,1}_-}$.

This completes the proof of estimate (\ref{2.3}).

The remaining cases $\pm_1 = -$ and $\pm_2 = \pm$ in (\ref{2.1}) and 
(\ref{2.1'}) can be treated in the same way. Using $\Pi_{\mp}(\eta) = 
\Pi_{\pm}(-\eta)$ we in fact get by (\ref{2.2})
\begin{eqnarray*}
\lefteqn{\left| \int \int \phi \langle \beta \Pi_-(D)\psi,\Pi_{\pm}(D)\psi' 
\rangle dx dt \right|} \\
& &= \left| \int \int \tilde{\phi} \langle \beta 
\Pi_-(\eta)\tilde{\psi}(\lambda,\eta),\Pi_{\pm}(\eta - 
\xi)\tilde{\psi}'(\lambda - \tau, \eta - \xi) \rangle d\lambda d\eta d\xi d\tau 
\right| \\
&&= \left| \int \int \tilde{\phi} \langle \Pi_{\pm}(\eta - \xi) \beta 
\Pi_+(-\eta)\tilde{\psi}(\lambda,\eta),\tilde{\psi}'(\lambda - \tau, \eta - 
\xi) 
\rangle d\lambda d\eta d\xi d\tau \right| \\
&& \le c \int \int \Theta_{\mp} |\tilde{\psi}(\lambda,\eta)| 
|\tilde{\psi}'(\lambda - \tau,\eta - \xi)| d\lambda d\eta 
|\tilde{\phi}(\tau,\xi)| 
d\tau d\xi \\
&& = I_{\mp} \, ,
\end{eqnarray*}
because by Lemma \ref{Lemma 3.1}
\begin{eqnarray*}
 \Pi_{\pm}(\eta - \xi) \beta \Pi_+(-\eta) & = & \sigma_{+,\pm}(-\eta,\eta - 
\xi) \\ 
& = & O(\angle(-\eta,\pm(\eta-\xi)) = O(\angle(\eta,\mp(\eta-\xi)) = 
O(\Theta_{\mp}) \, ,
\end{eqnarray*}
which can be handled like $I_{\pm}$ above, namely as follows. Our aim is to 
show
$$ I_{\pm} \le c \|\psi\|_{X^{0,\frac{1}{3},1}_-}  
\|\psi'\|_{X^{0,\frac{1}{3},1}_{\mp}} 
\|\phi\|_{X^{\frac{1}{2},\frac{1}{3},1}_{\pm 1}} \, . $$
This can be handled in the same way as before, provided the following Lemma 
holds.
\begin{lemma}
Denoting
 $$A_{\pm 1}= \tau \pm_1 |\xi| \, , \, B_- = \lambda - |\eta| \, , \, C_{\pm} = 
\lambda - \tau \pm|\eta-\xi| $$
 we have
 $$ \rho_{\pm} \le |A_{\pm 1}| + |B_-| + |C_{\mp}| $$
 where 
 $$ \rho_+ = |\xi| - ||\eta|-|\eta-\xi|| \, , \, \rho_- = |\eta| + |\eta-\xi| - 
|\xi| \, . $$
 \end{lemma} 
 {\bf Proof:}
 \begin{eqnarray*}
 \rho_+ & \le & |\xi| \pm |\eta| \mp |\eta - \xi| = |\xi| \pm \tau \mp\lambda 
\pm|\eta| \mp \tau \pm\lambda \mp|\eta-\xi| \\
 & \le & ||\xi|\pm\tau| + |\lambda-|\eta|| + |\lambda-\tau-|\eta-\xi|| = 
|A_{\pm}| + |B_-| + |C_-| 
 \end{eqnarray*}
and
\begin{eqnarray*}
\rho_- & = & |\eta| + |\eta-\xi| - |\xi| \\
& = & -\lambda + |\eta| + \lambda - \tau + |\eta-\xi|+\tau - |\xi| \le |B_-| + 
|C_+| + |A_-|
\end{eqnarray*}
as well as for $\tau \le 0$:
\begin{eqnarray*}
\rho_- & = & |\eta|+|\eta-\xi|-|\xi| = 
|\eta|-\lambda+\lambda-\tau+|\eta-\xi|+\tau-|\xi| \\
& \le & ||\eta|-\lambda| + |\lambda-\tau+|\eta-\xi|| \\  & \le &  |B_-|+|C_+|
\end{eqnarray*}
and for $\tau \ge 0$
\begin{eqnarray*}
\rho_- & = & |\eta| + |\eta-\xi| - |\xi| = 
|\eta|-\lambda+\lambda-\tau+|\eta-\xi|+\tau-|\xi| \\
& \le & ||\eta|-\lambda| + |\lambda-\tau+|\eta-\xi|| + \tau + |\xi| \le |B_-| + 
|C_+| + |A_+| \, .
\end{eqnarray*}
This completes the proof of the Lemma and Prop. \ref{Prop. 2.3}.
\section{A bilinear Strichartz' type estimate}
The following bilinear refinement is crucial for the estimate of the term 
$I_1^-$. It follows from the following proposition, which can be found in  
\cite{AFS}.\\
Defining 
$$ [(f,g)_{HH\to L}]^{\widehat{}}(\xi) := \int_{\bf R^2} \chi_{\{|\xi| << 
|\eta| 
+ |\xi - \eta|\}} \widehat{f}(\eta) \widehat{g}(\xi - \eta) d\eta \, , $$
where $\chi_A$ is the characteristic function of the set $A$.
\begin{prop}
\label{Prop. 3.1}
(\cite{AFS1}, Theorem 6)\\
Let 
$$ u_{\pm}(t):= e^{\mp it|D|} f $$
and
$$ v_{\pm}(t):= e^{\mp it|D|} g \, . $$
Then we have
$$ \| |D|^{-s_3} (u_{\pm},v_{\pm})_{HH \to L} \|_{L^2_{xt}} \le c 
\|f\|_{\dot{H}^{s_1}} \|g\|_{\dot{H}^{s_2}} \, , $$
where $ s_1+s_2+s_3 = \frac{1}{2} $ , $ s_1,s_2 < \frac{5}{8} $ , $ s_1 + s_2 > 
0 $.
\end{prop}
Using Prop. \ref{Prop. 1.2} we get   
\begin{Cor}
Under the assumptions of Prop. \ref{Prop. 3.1} the following estimate holds
$$ \| |D|^{-s_3} (u,v)_{HH \to L} \|_{L^2_{xt}} \le c 
\|u\|_{X^{s_1,\frac{1}{2},1}_{\pm}} \|v\|_{X^{s_2,\frac{1}{2},1}_{\pm}} \, $$
where it is essential that the two signs on the right hand side are equal.  
\end{Cor}
The following consequence is exactly what we  need in order to control $I_1^-$ 
in a suitable way.
\begin{prop}
\label{Prop. 3.2}
$$\| \langle D \rangle^{-\frac{1}{2}} (u,v)_{HH \to L} \|_{L^2_{xt}} \le c 
\|u\|_{X^{\frac{1}{6},\frac{1}{3},1}_{\pm}} 
\|v\|_{X^{\frac{1}{6},\frac{1}{3},1}_{\pm}} \, . $$
\end{prop}
{\bf Proof:}
The previous corollary is applied with $s_1=s_2=s_3= \frac{1}{6}$ leading to
\begin{equation}
\label{3.1}
\| \langle D \rangle^{-\frac{1}{6}} (u,v)_{HH \to L} \|_{L^2_{xt}} \le c 
\|u\|_{X^{\frac{1}{6},\frac{1}{2},1}_{\pm}} 
\|v\|_{X^{\frac{1}{6},\frac{1}{2},1}_{\pm}} \, . 
\end{equation}
It is interpolated with the following estimate which follows from Sobolev and 
the estimate $ \|f\|_{L^4_t(L^2_x)} \le c \|f\|_{X^{0,\frac{1}{4}}_{\pm}}$, 
which is proven like (\ref{!}).
\begin{eqnarray}
\label{3.2}
\| \langle D \rangle^{-\frac{2}{3}}(u,v)_{HH \to L} \|_{L^2_{xt}}  \le  c 
\|uv\|_{L^2_t(L^{\frac{6}{5}}_x)}\le c \|u\|_{L^4_t(L^{\frac{12}{5}}_x)}  
\|v\|_{L^4_t(L^{\frac{12}{5}}_x)} \\ \nonumber
 \le  c \|u\|_{L^4_t(H^{\frac{1}{6}}_x)} \|v\|_{L^4_t(H^{\frac{1}{6}}_x)}
 \le c \|u\|_{X^{\frac{1}{6},\frac{1}{4}}_{\pm}} 
\|v\|_{X^{\frac{1}{6},\frac{1}{4}}_{\pm}}
\le c \|u\|_{X^{\frac{1}{6},\frac{1}{4},1}_{\pm}} 
\|v\|_{X^{\frac{1}{6},\frac{1}{4},1}_{\pm}}
\end{eqnarray}
Complex bilinear interpolation between (\ref{3.1}) and (\ref{3.2}) gives the 
result, using
$$(X^{\frac{1}{6},\frac{1}{2},1}_{\pm},X^{\frac{1}{6} 
,\frac{1}{4},1}_{\pm})_{[\frac{1
}
{3}]} = X^{\frac{1}{6},\frac{1}{3},1}_{\pm} \, . $$

\end{document}